\documentclass[12pt]{article}
\usepackage{graphicx}
\usepackage{amsmath}
\usepackage{amsfonts}
\usepackage{amssymb}

\newcommand{\n}{\noindent}

\newcommand{\real}{\mbox{I${\!}$R}}
\newcommand{\intZ}{\mbox{Z${\!\!}$Z}}

\begin{document}

\title{A formula for the Whitehead group of a three-dimensional
crystallographic group}
\author{A. Alves and P. Ontaneda}
\date{October 9, 2003}
\maketitle

\begin{abstract}
{We give an explicit  formula for Whitehead group of
a three-dimensional crystallographic group $\Gamma$ in terms of
the Whitehead groups of the virtually infinite cyclic subgroups of
$\Gamma$. }
 \end{abstract}

\section{ Introduction }

\hspace{0.5cm}The $K$-theory of integral group rings of crystallographic
groups has been
studied by Connolly and Ko\'zniewski \cite{CK}, Farrell and Hsiang \cite{FH},
Farrell and Jones \cite{FJ1}, 
L\"{u}ck and Stamm \cite{LS}, Quinn \cite{Q2}, Tsapogas \cite{T}. 
Also, in \cite{P}, Pearson gives
explicit computations of the lower algebraic $K$-theory of
two-dimensional crystallographic groups.\\

In this article we give a simple formula for Whitehead
group of a three-dimensional crystallographic group $\Gamma$ in
terms of the Whitehead groups of the virtually infinite cyclic subgroups
of $\Gamma$. 
Here is our main result:\\

\n{\bf Main Theorem:} \textit{ Let  $\Gamma$ be a
3-crystallographic group. Then}

$$Wh(\Gamma)=\bigoplus_{G\in
{VC}_{\infty}(\Gamma)}Wh(G).$$
\textit{Moreover,  the direct sum in the formula above is finite.}\\

In the theorem above $VC_{\infty}(\Gamma )$ is the set of conjugacy classes
of maximal virtually infinite cyclic subgroups of $\Gamma$.\\

 \n{\bf Corollary 1:} \textit{ Let  $\Gamma$ be a 3-crystallographic
group. Then $Wh(\Gamma)$ is infinitely generated if and only if
$\Gamma$ contains a maximal  virtually infinite cyclic subgroup
$G$ with $Wh(G)$ infinitely generated.}\\

The following example answers a question proposed by F. T. Farrell:
is there any 3-crystallographic group whose Whitehead group is 
infinitely generated?\\

\n{\bf Example:} Consider the 2-crystallographic group $Pmm$.
This group is generated by translations $\tau_{e_i}(p)=p+e_{i},\,
p\in \real^{2}$, $e_{1}=(1,0),\, e_{2}=(0,1)$ and reflections
about the $x$ and $y$ axis. Then $\intZ_{2}\times\intZ_{2}$ is a
subgroup of $Pmm$. Hence $\intZ_{2}\times\intZ_{2}\times\intZ$ is
a subgroup of the 3-crystallographic group
$\Gamma=Pmm\times\intZ$. In fact,
$\intZ_{2}\times\intZ_{2}\times\intZ$  is maximal because
$\intZ_{2}\times\intZ_{2}\times\intZ=\Gamma^{({\tiny{\real}}{e_{3}})}$
(see Lemma 5.1) and the action of $\intZ_{2}\times\intZ_{2}\times\intZ$
on $\real{e_{3}}$ is cocompact. Also by the
Bass-Heller-Swan formula 
$Wh(\intZ_{2}\times\intZ_{2}\times\intZ)=2\,
Nil_1(\intZ[\intZ_{2}\times\intZ_{2}])$, which is infinitely generated
(see Lemma 6.2). Therefore Corollary 1 implies that
$Wh(Pmm\times\intZ)$ is
infinitely generated.\\

\n{\bf Corollary 2:}  \textit{Let $\Gamma$ be a
2-crystallographic group. Then}
$$2\, Nil_{1}(\intZ[\Gamma])=2\left(\bigoplus_{F\in
F(\Gamma)}Nil_{1}(\intZ [F])\right).$$
\vspace{1cm}

In the corollary above $F(\Gamma )$ is the set of conjugacy classes
of maximal finite subgroups of $\Gamma$.\\

As in \cite{BFJP}, \cite{BJP}, \cite{LS}, \cite{P}, \cite{U}, 
the proof of our main theorem is accomplished using fundamental
results of Farrell and Jones \cite{FJ1}, ``the Farrell-Jones Isomosphism
Conjecture'', which holds for crystallograsphic groups. This
conjecture states that the algebraic $K$-theory of $\intZ[\Gamma]$ may be
computed from the algebraic $K$-theory of the virtually cyclic subgroups
of $\Gamma$ via an appropriate ``assembly map''.
Our approach is quite geometric and our main ingredient is the construction
of a somehow ``concrete" model for the universal $(\Gamma ,VC(\Gamma ))$ space,
for any $n$-crystallographic group.
This construction is a variation of Farrell and Jones construction given in
\cite{FJ1} and \cite{FJ3}. It is interesting to note that some the results for 
2-crystallographic groups in \cite{P} can be obtained in a simpler fashion 
using our geometric methods.\\

In this article we also classify, modulo isomorphism, the virtually infinite cyclic
subgroups of a 3-crystallographic group.\\

This paper is organized as follows. In the Section 1, we 
introduce some definitions and
propositions and state Farrell-Jones Isomorphism Conjecture.
In the Section 2, we construct a model for the
universal $(\Gamma,VC(\Gamma))$-space, where $\Gamma$ is a
$n$-crystallographic group. In the Section 3, we classify,
modulo isomorphism, the virtually  infinite cyclic subgroups of a
3-crystallographic group. In the Section 4, we  calculate the
isotropy  groups  of open $n$-cells of the universal
$(\Gamma,VC(\Gamma))$-space. In the Section 5, we calculate the
${\cal E} ^2_{i,j}$ terms of the spectral sequence, necessary to
apply the Farrell-Jones  Isomorphism Conjecture. Finally, in Section 6,
we proof the results mentioned above.\\

We wish to thank Prof. F. T.
Farrell for suggesting the problem and for the helpful information about
$Wh(A_4\times\intZ_2)$. We are also grateful to Profs. D. Juan-Pineda and S.
Prassidis for their helpful information about $K$-theory of finite groups.

\section{ Preliminaries }

\hspace{0.6cm}First, let us recall some definitions and fix some notation.\\

A group $\Gamma < O(n)\rtimes\real^n$ is
n-crystallographic group if $\Gamma$ acts properly
discontinuously,
cocompactly and by isometries on $\real^n$.\\

Let  $\Gamma$ be a group and $G<\Gamma$, i.e., $G$ is a subgroup of
$\Gamma$. $G$ is a
maximal finite subgroup of $\Gamma$, if $G<G'<\Gamma$, with $G'$
finite, implies $G=G'$. Analogously, we define a maximal virtually
infinite cyclic subgroup of $\Gamma$.\\

Let $X$ be a space equipped with an action of a group $\Gamma$. Let
$B\subset X$. Then we write\\

 $\Gamma^B=\{\gamma\in\Gamma;
\gamma b=b$ for all  $b\in B \}$.\\

 $\Gamma^{(B)}=\{\gamma\in\Gamma; \gamma B=B\}$.\\

 If $G<\Gamma$,
$X^G=\{x\in X; gx=x$ for all $g\in G\}$.\\

 If $Y$ is another  $\Gamma$-space, we write $X{\cong}_{\Gamma}Y$
if $X$ is
$\Gamma$-homeomorphic to $Y$.\\

A group $G$ is virtually cyclic if it is finite or
contains an infinite cyclic subgroup of finite index.\\

For a group $\Gamma$, we write $F(\Gamma)$ for the set of
conjugacy classes of maximal finite subgroups of  $\Gamma$, and
$VC_{\infty}(\Gamma)$ the set of conjugacy classes of maximal
virtually infinite cyclic subgroups of  $\Gamma$.
$VC(\Gamma)$ will denote the family of
virtually cyclic  subgroups of $\Gamma$.\\

We will use the following result.\\

\n{\bf Proposition 1.1}  \textit{Let  $G$ be
a subgroup of a $n$-crystallographic group $\Gamma$. Then $G$ virtually infinite 
cyclic if and only if there is at least one line $l$ in $\real^n$
left invariant by $G$,  and the action of $G$ on $l$ is cocompact.}\\

{\bf Proof.} This result is proven, with greater generality, in \cite{FJ1}
p. 267. The proposition can also be shown directly, using some simple geometric
arguments and the algebraic result of \cite{SW} p. 178-179.\\

Now, we shall state Farrell-Jones Isomorphism Conjecture (see \cite{FJ1}).
Let $X$ be a connected CW-complex and $\Gamma=\pi_1(X)$. Let $A$ denote a universal
$(\Gamma,VC(\Gamma))$-space and $\tilde X$ be the universal covering space
of $X$. Let $\Gamma\times \tilde X\times A\longrightarrow  \tilde X\times
A$ be the diagonal action. Then $\rho:{\cal E}(X) \rightarrow{\cal B}(X)$ and
$f:{\cal E}(X)\rightarrow X$ are defined to be the quotient of the
standard projections $\tilde X\times
A \rightarrow A$ and  $\tilde X\times
A \rightarrow \tilde X$ under the $\Gamma$-actions. Let $\zeta_*()$ denote
any of the
$\Omega$-spectra-valued functors of \cite{FJ1} p.251-253. Then Farrell and Jones conjecture
that $${\mathbb H}_{*}({\cal B}(X),{\zeta}_{*}(\rho))
\stackrel{{\zeta}_{*}(f)\circ A_{*}}
{\longrightarrow}
{\zeta}_{*}(X)$$
is an equivalence of the $\Omega$-spectra, where $A_*$ is an ``assembly
map''
for the simplicial stratified fibration $\rho:{\cal E}(X) \rightarrow{\cal
B}(X)$ (see \cite{FJ1} p.257) and ${\zeta}_{*}(f)$ is the image of $f:{\cal
E}(X)\rightarrow X$ under the functor ${\zeta}_{*}()$.\\

Let ${\cal P}_*()$ denote the functor that maps $X$ to the $\Omega$-spectrum 
of stable topological pseudoisotopies on $X$. The following result was proven
in \cite{FJ1}:\\

\n{\bf Farrell-Jones Isomorphism Proposition for ${{\cal P}_*()}$
1.2:} \textit{The above conjecture is true for the functor ${\cal
P_*()}$ on the space $X$ provided that there exits a simply
connected symmetric Riemannian manifold $M$ with non positive
sectional curvature everywhere such that $M$ admits a properly
discontinuous cocompact group action of $\Gamma=\pi_1(X)$ by
isometries of $M$.}\\

Follows that the Farrell-Jones Isomorphism Conjecture holds for crystallographic
groups.\\

The relationship between ${\cal P_*()}$ and lower algebraic $K$-theory is
given by\\

\n{\bf Proposition 1.3:} \cite{AH}    \textit{$$ \pi_{j}({\cal
P}_{*}(X))=\left\{\begin{array}{rcl}
K_{j+2}(\intZ[\pi_{1}(X)]), & \mbox{\it if} & j\leq -3\\
{\tilde K}_{0}(\intZ[\pi_{1}(X)]), & \mbox{\it if} & j=-2\\
Wh(\pi_{1}(X)), & \mbox{\it if} & j=-1 .
\end{array}
\right.
$$ } \\

Finally, the following result aids in the calculation of the lower $K$-theory
of $\Gamma$.\\

\n{\bf Proposition 1.4: }  (see \cite{Q1})\textit{ Let $f: E \rightarrow X$
be a simplicially stratified fibration. Then there is a spectral sequence
with ${\cal E} ^{2}_{i,j}=H_{i}(X,\pi_{j}\zeta_{*}(f))$ which abuts to
$H_{i+j}(X,\zeta_{*}(f))$.} \\

\section{ A model for the universal $(\Gamma,VC(\Gamma))$-space}

\hspace{0.5cm} In this section we construct a model for the universal
$(\Gamma,VC(\Gamma))$-space, where  $\Gamma$ is a crystallographic
group and $VC(\Gamma)$ is the set of virtually cyclic subgroups of
$\Gamma$.\\

Let $\Gamma$ be a n-crystallographic group. Then $\Gamma$ fits in
the short exact sequence of groups $$0\rightarrow
T\stackrel{i}\rightarrow \Gamma\stackrel{p}\rightarrow
F\rightarrow 1$$ with $T$ isomorphic to $\intZ^{n}$ and $F$ a finite
subgroup of $O({n})$.
Every $\gamma\in\Gamma$ acts on $\real^n$ as:\\

 $\gamma x=Jx+b$, where $J=p\gamma\in O(n)\,\, b\in\real^n.$\\
 
\n for all $x\in \real^n$. Note that\\

 $\gamma^{-1} x=J^{-1}x-J^{-1}b.$\\

The elements in $T$ act on $\real^n$ by
translations. For $h\in \real^n$ we write $\tau_h$ for the translation by
$h$, i.e., $$\tau_h:\real^n\rightarrow \real^n, \, \, \tau_h x=x+h.$$
Let $H=\{h;\, \tau_h\in T\}$. Then  $H<\real^{n}$ is
an additive subgroup of
$\real^n$. \\

Remark that, if $h\in H$ and $\gamma\in \Gamma$, then $\,
(p\gamma)h\in H$. To see this, first notice that $\gamma
\tau_{h}{\gamma}^{-1}\in\Gamma$.  Then $$(\gamma
\tau_h\gamma^{-1})x=x+J h=\tau_{J h}x\Rightarrow Jh=(p\gamma)h\in H.$$

Define $C^{n}:=\real^{n}/\sim$, where $x\sim y\Leftrightarrow
x=\pm y$. Hence $C^{n}$ is homeomorphic to  $C(\real P^{n-1})$,
the open cone over the real projective space $\real P^{n-1}$. The
\textit{vertex} of $C^{n}$ is the point $[\, 0\, ]\in C^{n}$. A
subset of the form $\{[\, \lambda x\, ], \lambda\in \real\}$, $x\neq
0$ is
called a \textit{ray} of $C^{n}$.\\

 Note that every element $[x]\in C^{n}$
\vspace{0,1cm}
has a well defined ``norm'' $||\, [x]\, ||=||x||=||-x\, ||$, that
satisfies
$||\, [\lambda x]\, ||=|\lambda|\, ||x||,\lambda\in\real.$\\

Note also that $\Gamma$ acts on  $C^{n}$ via $p$: $$\gamma
[x]=[(p\gamma) x].$$
Consequently $\Gamma$ acts on
$C^{n}\times\real^{n}$ diagonally: $$\gamma ([x],y)=(\gamma[x],\gamma
y)=([(p\gamma )x],\gamma y).$$

For each  $h\in H<\real^{n}, h\neq 0$, write ($\, \parallel$ means ``parallel")\\

 $l_{h}={\real}{h}=\{th, t\in\real\}$.\\

 ${\cal L}_{h}=\{l$, line in $\real^{n};l\parallel l_{h}\}$.\\

 ${\cal L}=\bigcup_{h\in H}{\cal L}_{h}$.\\

 $\Lambda=\{{\cal L}_{h},h\in H\}$.\\

\n Note that
${\cal L}_{h}$,  with the
quotient topology, is
homeomorphic to $\real^{n-1}$.
The action of $\Gamma$ on $C^n$ induce actions of $\Gamma$ on\\

\n{\bf (1)}  ${\cal L}$, because $l\parallel l_{h}$ $\Rightarrow$
$\gamma
l\parallel l_{(p\gamma) h}$, for $h\in H,\gamma\in \Gamma$.\\

\n{\bf (2)} $\Lambda$, by $\gamma{\cal L}_{h}={\cal L}_{(p\gamma)h}$.\\

\n Note that
the previous remark implies that the action is well-defined. Also, if
$\gamma$ is a translation, $$\gamma{\cal L}_{h}={\cal
L}_{h},\forall h\in H.$$\\

Since $p\gamma\in F$ and $F$ is finite, we have that the  orbit
$\Gamma{\cal L}_{h}$ of ${\cal L}_{h}$ is finite, $\forall h\in
H$. Enumerate these orbits: $\Lambda_{1}, \Lambda_{2},
\Lambda_{3},...$
then:\\

\n {\bf (1)} Each $\Lambda_{k}$ is finite.\\

\n{\bf (2)} $\Gamma\Lambda_{k}=\Lambda_{k}$.\\

\n{\bf (3)} $\bigcup_{k}\Lambda_{k}=\Lambda$.\\

For each $l\in {\cal L}$ define $c(l)\in C^{n}$, the {\it height} of $l$,
in the following
way: If $l\in{\cal L}$, then $l\parallel l_{h}$ for some $h\in H$.
Hence there is $k$ such that ${\cal L}_{h}\in \Lambda_{k}$.  Define
\vspace{0,1cm}
$c(l)={\big [}\frac{kh}{\parallel h\parallel}{\big ]}\in C^n$.
Note that $c(l)$ is well defined. Note also that
$\gamma c(l)=c(\gamma l)$ and that $c$ is constant on each ${\cal
L}_{h}$.\\

For $l\in{\cal L}$, write $\bar l=\{c(l)\}\times l\subset
C^n\times\real^n$. We think of $\bar l$ as the line
$l\subset\real^n$ lifted to a height $c(l)\in C^n$. Define
$\bar{\cal L}=\{\bar l,l\in{\cal L}\}$, and $\bar{\cal
L}_{h}=\{\bar l,l\in{\cal L}_{h}\}$. Note that $\bar{\cal
L}=\bigcup_{h\in H}\bar{\cal L}_{h}$. Then $\bar{\cal L}$ is a set
\vspace{0,1cm}
of disjoint ``lifted lines" in $C^{n}\times\real^{n}$. 

%That is, if ${\bar
%l}_{1}, {\bar l}_{2} \in \bar{\cal L}$, then either ${\bar
%l}_{1}\cap{\bar l}_{2}=\emptyset$ or ${\bar l}_{1}={\bar l}_{2}$.\\

Finally, define $A=C^{n}\times\real^{n}/\cong$, where $x\cong
\vspace{0,1cm}
y\Leftrightarrow x=y$ or $x,y\in \bar l$, for some ${\bar
l}\in{\bar{\cal L}}$, i.e., $A$ is obtained from
$C^{n}\times\real^{n}$ by identifying each line in $\bar {\cal L}$
to a point.\\

Let $Y=\bigcup_{\bar l\in\bar{\cal L}}\bar l \subset
C^n\times\real^n$. Then $\Gamma Y=Y$ and
\vspace{0,1cm}
$\Gamma((C^{n}\times\real^{n})-Y)=(C^{n}\times\real^{n})-Y$. Hence the
action of $\Gamma$ on $C^{n}\times\real^{n}$ induces an action of
$\Gamma$ on $A$.
 Let $\pi:C^{n}\times\real^{n}\rightarrow A$
be the collapsing map, and write $Z=\pi Y$. Then $\Gamma Z=Z$ and
$\Gamma (A-Z)=A-Z$. Note that
$\pi:(C^{n}\times\real^{n})-Y\rightarrow A-Z$
 is a  $\Gamma$-equivariant homeomorphism.
Since $\bar{\cal L}=\bigcup_{h\in H}\bar{\cal L}_{h}$, we have
$Y=\bigcup_{h\in H}(\bigcup_{\bar l\in \bar{\cal L}_{h}}\bar l)$.
\vspace{0,1cm}
Consequently $Z=\pi Y=\bigcup_{h\in H}\pi Y_{h}=\bigcup_{h\in
H}Z_{h}$, where $Y_{h}=\bigcup_{\bar l\in \bar{\cal L}_{h}}\bar l$
and $Z_{h}=\pi Y_{h}$.\\

\centerline{\large{\bf The action of $\Gamma$ on $Z$}}

 \vspace{0,5cm}

Since the connected components of $Z$ are the $Z_h$'s, we have
that either $\gamma{Z}_h={Z}_h$ or
$\gamma{Z}_h\cap{Z}_h=\emptyset$, for $\gamma\in\Gamma$. Consider
$\Gamma^{({Z}_h)}=\{\gamma\in\Gamma; \gamma{Z}_h={Z}_h\}$. Note
that $\Gamma^{({Z}_h)}=\Gamma^{({\bar{\cal L}}_h)}=\Gamma^{({\cal
L}_h)}$. We want to study the action of $\Gamma^{({Z}_h)}$ on
${Z}_h$, or equivalently, the action of $\Gamma^{({\cal L}_h)}$ on
${\cal L}_h$. For this we need an explicit homeomorphism between
${\cal L}_h$ and $\real^{n-1}$. Define
$\alpha_h:=(l_h)^{\perp}\subset\real^n$, i.e., $\alpha_h$ is the
$(n-1)$-space orthogonal to the line spanned by $h$. Then for
every $l\in{\cal L}_h$, we define $p_l$ to be the unique point
such that $\{p_l\}:=l\cap\alpha_h$. Write $\varphi(l)=p_l$. It is
easy to verify that $\varphi$ is a homeomorphism between ${\cal
L}_h$
and $\alpha_h\cong\real^{n-1}$.
Let $\gamma\in\Gamma$. Then there are $J\in O(n), a\in\real^n$
such that $\gamma x=J x+a$, for all $x\in\real^n$. We write
$\gamma= (J,a)$. Then $\gamma\in
\Gamma^{({\cal L}_h)}$ if and only if $J h=\pm h$.
This implies that $J l_h=l_h$ and $J\alpha_h=\alpha_h$. Hence,
$J|_{\alpha_h}:\alpha_h\rightarrow \alpha_h$ and $J|_{\alpha_h}\in
O({\alpha_h})$ (${\alpha_h}$
with the scalar product of $\real^n$).
Write $a=a_0+\lambda h$, with $a_0\in\alpha_h$, i.e.,
$a_0={Proj}_{\alpha_h}a$. Define
$\bar\gamma:{\alpha_h}\rightarrow{\alpha_h}$,
$\bar\gamma:=\varphi\gamma\varphi^{-1}$. Hence we have the
commutative diagram:

$$\begin{array}{ccc}
{\cal L}_h & \stackrel{\gamma}{\longrightarrow} &
 {\cal L}_h \\
 {\downarrow}{\varphi} & & {\downarrow}{\varphi} \\
 \alpha_h  & \stackrel{\bar\gamma}{\longrightarrow} & \alpha_h \\

\end{array}
$$\\

A simple calculation shows that $\bar\gamma=(J|_{\alpha_h},a_0)$,
hence $\Gamma^{({\cal L}_h)}$ acts by isometries on
${{\cal L}_h}$, where we consider ${\cal L}_h$ as a vector space
with scalar product obtained by identifying ${\cal L}_h$ with
$\alpha_h\subset\real^n$.
Moreover, it is not difficult to show that this action is
crystallographic, i.e., properly discontinuous and cocompact 
(this is proved, with more generality, in \cite{FJ1} p.267).\\

\vspace{1cm}

 \centerline{\large{\bf The triangulation of $A$}}

 \vspace{1cm}

 $C^{n}$ and $\real^{n}$ are $PL$
$\Gamma$-spaces, hence $C^{n}\times\real^{n}$ is a $PL$
$\Gamma$-space. Let ${\cal T}_C$ be a $\Gamma$-equivariant
triangulation of $C^{n}$ and ${\cal T}_R$ a $\Gamma$-equivariant
triangulation of $\real^{n}$. Let $\tilde{\cal T}$ be the cell
structure on $C^n\times\real^{n}$ with product cells
$\sigma_C\times \sigma_R, \sigma_C\in {\cal T}_C, \sigma_R\in{\cal
T}_R$. Then $\tilde{\cal T}$ is $\Gamma$-equivariant. We can
suppose that each $Y_{h}$ is a subcomplex of
$C^{n}\times\real^{n}$ (because the projection of each $Y_{h}$ into 
$C^{n}$ is a point). Since the triangulations  ${\cal T}_C$ and 
${\cal T}_R$ are $\Gamma$-equivariant, we have that
$\tilde {\cal T}$ induces a $\Gamma$-equivariant
cell structure ${\cal T}$ on $A$
with the following properties:\\

\n{\bf (1)} Each $Z_{h}$ is a subcomplex of $A$.\\

\n{\bf (2)} If $\sigma$ is an open k-cell in $A-Z\cong_{\Gamma}
(C^n\times\real^{n})-Y$, then $\sigma=\sigma_C^i\times\sigma_R^j$,
$k=i+j$, where $\sigma_C^i$ and $\sigma_R^j$ are open simplices
in $C^n$ and $\real^n$, respectively.\\

We say that a 1-cell $\sigma^1_C$ in $C^n$ is a  \textit{ray cell} if
$\sigma^1_C$ is contained in a ray of $C^n$.\\

\n{\bf Theorem 2.1:}  \textit{$\Gamma\times A\longrightarrow A$ is an
universal
$(\Gamma,VC(\Gamma))$-space.}\\

\n{\bf Proof:} We must  verify the following properties:\\

\n{\bf (1)} $A$ can be equipped with a cell structure $K$ such that
$\Gamma\times A\longrightarrow A$ is a cellular action. Moreover, for
$\gamma\in\Gamma$ and $\sigma \in K$ we have that $\gamma(\sigma)=\sigma$,
implies $\gamma|_{\sigma}$=inclusion.\\

\n{\bf (2)} For any $p\in A$ we have that $\Gamma^p\in VC(\Gamma)$.\\

\n{\bf (3)} For each subgroup $G\in VC(\Gamma)$ we have that $A^G$ is a
nonempty contractible subcomplex of $K$.\\

(1) follows from the discussion above, taking $K=\cal{T}$.\\
(To get that $\gamma(\sigma)=\sigma$,
implies $\gamma|_{\sigma}$=inclusion, we may have to subdivide
$\cal{T}_C$ and $\cal{T}_R$.)
To verify (2) we first note that if $p\in A-Z$ then  $\Gamma^p$ is
a finite group. On the other hand, if $p\in Z$, then $p\in Z_h$,
 for some $h\in H$. Hence  $\Gamma^p=\{\gamma\in\Gamma: \gamma(q+\real
h)=q+\real
h\}$, where $\pi^{-1}(p)=q+\real h$. Thus $\Gamma^p$ acts properly
discontinuously cocompactly by isometries on the line $q+\real h$. Consequently, 
by prop. 1.1, $\Gamma^p\in VC(\Gamma)$.\\

Finally  we verify (3). First consider the case where $G\in
VC(\Gamma)$ contains an element $g\in G$ of infinite order, i.e.,
$G$ contains translations. Then
$A^G\cap((C^n\times\real^n)-Y)=\emptyset$. Hence $A^G\subset Z$.
Moreover $A^G\subset Z_h$, where $h\in H$ is such that $\tau_h\in
G$. Let $l=a+\real h$ be the line such that $Gl=l$, i.e.,
$G<\Gamma^{(l)}$ (see prop. 1.1). Also $[\, l\, ]\in Z_h\subset Z\subset A$. 
Since $Gl=l$, we have that $G[\, l\, ]=[\, l\, ]$, that is $[\, l\, ]\in A^G$.
Follows that $A^G\neq\emptyset$.\\

Note that $G\subset\Gamma^{(Z_h)}$ and $A^G\subset Z_h$ is
precisely the fixed point set of the $G$ action on $Z_h$. We know
that the action of $\Gamma^{(Z_h)}$ on $Z_h$ is linear (in fact, it is
crystallographic). Hence $A^G=\bigcap_{g\in G}A^{\{g\}}$, and each
$A^{\{g\}}$ is a vector space. Therefore $A^G$ is contractible.\\

Consider now the case where $F\in VC(\Gamma)$ is finite.\\ 

Suppose first that $F$ is trivial, i.e. $A^\Gamma$=$A$. We prove that
$A$ is contractible.
(To prove this we cannot just use the fact that the collapsing map
$\pi : C^{n}\times \real^{n}\rightarrow A$ is cell-like, because $\pi$ is not
proper.)\\

For $0\neq x\in\real^n$, let $\alpha_x$ be the plane orthogonal to $\real x$ in
$\real^n$, and for $v\in \real^n$ define $proj_x v$ to be the orthogonal projection
of $v$ in $\alpha_x$.
We define a homotopy $h_t$ on $C^n \times\real^n$ that moves a point $([x],v)\in
C^n \times\real^n$, $x$ away from $0$, toward $([x],proj_x v)$:
$$h_{t} ([x],v)=\bigl( \,\,[x]\,\, ,\,\, (\delta(||x||)t)proj_x v\, +\, (1-\delta(||x||)t)v
 \bigr)$$

\n where $\delta :[0,\infty )\rightarrow [0,1]$ is a continuous function such that
$\delta (t)=0$ for $0\leq t\leq 1/3$, $\delta (t)=1$ for $2/3\leq t$ and 
$0<\delta (t)<1$ for $1/3<t< 2/3$.\\

Then $h_0 $ is the identity and $h_1 (C^n \times\real^n ) - (B_{2/3} \times\real^n )=
\{([x],v)\,\, ,\,\, v\in\alpha_x\, ,\,\, \,\,\,2/3\leq ||x||\}$. Here $B_{2/3}=\{ x\in C^n \,\, ,\,\, ||x||<2/3\}$.\\

It is straightforward to prove that $h_t : C^n \times\real^n
\rightarrow C^n \times\real^n$ defines a homotopy $H_t :A\rightarrow A$, with 
$H_t \pi = \pi h_t$, where $\pi : C^{n}\times \real^{n}\rightarrow A$ is the collapsing map.\\

Since the collapsing occurs away from $B_2/3 \times\real^n$, we have that 
$\pi |_{h_{1}(  C^{n}\times \real^{n})}$ is an embedding. Hence
$$A=H_0 A\sim H_1 A=H_1 \pi ( C^{n}\times \real^{n})=\pi h_1 
(C^{n}\times \real^{n})\cong_{TOP} h_1 (C^{n}\times \real^{n})$$
$$\sim h_0 (C^{n}\times \real^{n})
=C^{n}\times \real^{n}\sim pt$$
($\sim$ means  ``homotopy equivalent" and $\cong_{TOP}$ means ``homeomorphic".)
This proves that $A$ is contractible.\\

Suppose now that $F$ is any finite subgroup of $\Gamma$. Hence
there is $x_0\in\real^n$ such that $Fx_0=x_0$.
Thus $([0],x_0)\in A^F$.\\

We prove that $A^F$ is contractible. Let
$B=\pi^{-1}(A^F)$. Without loss of generality,
we can suppose that $x_0=0\in\real^n$, that is, $F0=0$.
Then, every $f\in F$ is linear. 
Define a subspace $E_0$ and a set $E_1$ in the following way:\\

\n{\bf (1)}  $E_0 =\{ y\in\real^n\, ;\,\, fy=y$ for all $f\in F\}$.\\

\n{\bf (2)} $E_1=\{y\in\real^n; fy=\pm y$ for all $f\in F$ and
 there is $f\in F$ such that $fy=-y$\}.\\

Note that $E_1$ is a cone (i.e. $y\in E_1 $ implies $\lambda y\in E_1$, $\lambda\in\real$).
Note also that $E_0\cap E_1 =\{ 0\}$. Moreover, if $y\in E_1$ then $y\perp E_0$.\\

 Let $PE_i\subset \real P^{n-1}$ be the projectivization of $E_i\subset\real^n$
and $C(PE_i)\subset C(\real P^{n-1})= C^n$ the  cone of $PE_i$.\\

\n{\bf Claim 1:} \textit{$B=\big\{[C(PE_0)\cup C(PE_1)]\times
E_0\big\}\cup\big\{\bigcup_{l\in{\cal L}, Fl=l}\bar l\, \big\}$.}\\

\n{\bf Proof:} Consider $B_1=B\cap [(C^n\times\real^n)-Y]$ and
$B_2=B\cap Y$. Then $B=B_1\coprod B_2$. Note that if
$([x],y)\in [C(PE_0)\cup C(PE_1)]\times E_0$, then $Fy=y$ and
$Fx=\pm x$. Hence $([x],y)\in B$. Note also that $B_2=B\cap
Y=\bigcup_{l\in{\cal L}, Fl=l}\bar l\subset B$. Therefore
$$\big\{[C(PE_0)\cup C(PE_1)]\times
E_0\big\}\cup\big\{\bigcup_{l\in{\cal L}, Fl=l}\bar l\, \big\}\subset B.$$
Conversely if $([x],y)\in B_1$, then $Fy=y$ and $Fx=\pm x$. Thus
$y\in E_0$ and $x\in E_0$ or $E_1$. Hence $[x]\in C(PE_0)\cup
C(PE_1)$. Therefore $([x],y)\in [C(PE_0)\cup
C(PE_1)]\times E_0$. This proves the claim.\\

\n{\bf Claim 2:} \textit{ If $Fl=l$, then $l\cap E_0\neq\emptyset$.
Moreover either $l\subset
E_0$ or $l\subset E_1+y$, $\{y\}=l\cap E_0$.}\\

\n{\bf Proof:} Since $Fl=l$ and $F$ is finite, there is $y\in
l$, such that $Fy=y$. Hence $y\in E_0$, and follows that $l\cap
E_0=\{ y\}\neq\emptyset$.\\

 Suppose $l\not\subset E_0$. We know that  $l\parallel l_h$, for some
$h\in H$.
Then $y\pm h\in l$, where $\{y\}=l\cap E_0$. Since $fy=y$ for all
$f\in F$,  we have that $f(y+h)=y\pm
h$, for all $f\in F$. Follows that $f(h)=\pm h$, for all $f\in F$.
Since there is a $f\in F$ such that $f(y+h)=y-h$, follows
that $fh=-h$ for some $f\in F$. Hence $h\in E_1$. Therefore $l\subset
E_1+y$. This proves the claim.\\

If $l\subset E_0$, then $l\parallel l_h=\real h$, $h\in E_0$. Thus
$c(l)\in C(PE_0)$. Consequently, if $a\in l$ we have $(c(l),a)\in C(PE_0)\times
E_0$. Hence $\bar l\subset C(PE_0)\times E_0$. On the other hand,
 if $l\subset E_1+y$, $l\parallel l_h$, $h\in E_1$. Thus $c(l)\in
 C(PE_1)$. If $a\in l$, then $(c(l),a)\in [C(PE_0)\cup C(PE_1)]\times
E_0$ if and only if $a\in E_0$, i.e., $a=y$. Therefore $\bar
l\cap[C(PE_0)\cup C(PE_1)]\times E_0=\{point\}$.\\

Follows from the discussion above that collapsing every line $\bar l$,
with $l$ not in $E_0$, to its (unique) point of intersection with
$[C(PE_0)\cup C(PE_1)]\times E_0$, we obtain a homeomorphism of $A$ onto
the space $\pi([C(PE_0)\cup C(PE_1)]\times E_0)$. Moreover, 
since there is no collapsing among points in $C(PE_1) \times E_0$, deforming
$\pi ( C(PE_1)\times E_0)$ to $\{ [0]\}\times E_0$
we see that $A$ is homotopy equivalent to $\pi ( C(PE_0)\times E_0)$,
which is obtained from $C(PE_0)\times
E_0$ by collapsing to a point every line $\bar l \subset C(PE_0)\times
E_0$, with $l\subset E_0$, $l\in\cal L$. But an argument similar to the one used before 
(to prove that $A=\pi (C^n \times \real^n)$ is contractible) shows that
$\pi ( C(PE_0)\times E_0)$ is contractible. This
proves  property (3). Therefore $\Gamma\times A\longrightarrow
A$ is an universal $(\Gamma,VC(\Gamma))$-space. This proves the
theorem.

\vspace{1cm}

\section{ Virtually infinite cyclic subgroups in dimension three}

\hspace{0.5cm} In this section we calculate, modulo isomorphism, the
virtually
infinite
cyclic subgroups of a 3-crystallographic group.\\

In what follows, $\intZ_i$ will denote the
cyclic group of order $i$, $D_i$  the dihedral group of order
$2i$, $S_{n}$  the  permutation  group of order $n!$ and $A_n$ will 
denote the alternating group of order $n!/2$.\\

We will use the following facts
about a 3-crystallographic group $\Gamma$ and a finite subgroup $F$ of
$\Gamma$.\\

\n{\bf (a)} If $\gamma\in \Gamma $ has finite order, then $1\leq
|\gamma|\leq 6$ and $|\gamma |\neq 5$. This is  due to the crystallographic restriction
(see \cite{S}, p.32).\\

\n{\bf (b)} If $F$ is finite, then $F$ is trivial or isomorphic to
one of the following groups: $\intZ_i, D_i, \intZ_i\times\intZ_2,
D_i\times\intZ_2, A_4, S_4,
A_4\times\intZ_2, S_4\times\intZ_2,\, i=2,3,4,6$ (see \cite{S}, p.49).\\

\n{\bf (c)} If $F$ is isomorphic to one of the groups $A_4, S_4,
A_4\times\intZ_2, S_4\times\intZ_2$, then $Fl\neq l$ for all lines
in $\real^3$, i.e., $F$ does not leave any line in $\real^3$ invariant. In
particular $(\real^3)^F=\{x\in \real^3;\gamma x=x$,
for all $\gamma\in F\}=\{point\}$ (see \cite{S}, p.48).\\

\n{\bf (d)} If $F$ fixes a point $x_0$ and leaves invariant a line
$l$ that contains $x_0$ (or equivalently, a plane $\alpha$ that
contains $x_0$) then, by (c) and (b) above, follows that $F$ is trivial or
isomorphic to one of the groups $\intZ_i, D_i,
\intZ_i\times\intZ_2, D_i\times\intZ_2$, $
i=2,3,4,6.$\\

Now, let $G$ be a virtually infinite cyclic subgroup of a
3-crystallographic group $\Gamma$. Recall that $G$ leaves
invariant a line $l$ in $\real^3$, i.e.,
$G<\Gamma^{(l)}$, and $G$
acts cocompactly on  $l$ (see prop. 1.1). Hence we have two possibilities:\\

\n{\bf (1)}\, $G$ acts by translations and by reflections on $l$,
i.e., $G$ has
a dihedral action on $l$.\\

\n{\bf (2)}\, $G$ acts just by translations on $l$.\\

In the first case, we have a surjection $G\rightarrow
D_{\infty}$, and in the second case we have a
surjection $G\rightarrow\intZ$. In any case, we have

$$0\rightarrow F\hookrightarrow G\stackrel{\rho}\rightarrow H\rightarrow
0,$$
\n where $H$ is isomorphic to $D_{\infty}$ or $\intZ$,
$\rho(g)=g|_l:l\rightarrow l$, and $F$ is
the subgroup of all elements $g\in G$, that fix $l$ pointwise,
i.e., $g|_l=1_l$. Hence $F$ acts on a plane $\alpha$ orthogonal to
$l$ (choose any plane orthogonal to $l$). By (a) above we have that 
$F$ is trivial or isomorphic to $\intZ_i$ or $D_i$, $i=2,3,4,6$. 
We analyze the two cases separately.\\
\vspace{1cm}

\centerline{\large{\bf $H\cong D_{\infty}$}}

\vspace{0,5cm}

Then $G$ has a dihedral action on $l$. Write
$D_{\infty}=\intZ_2^a*\intZ_2^b$, where $a$ and $b$ are
reflections about some points $p_a$, $p_b\in l$, $p_a\neq p_b$. 
Thus $G=G^a*_{F}G^b$, where $G^a=\rho^{-1}(\intZ_2^a)$ and 
$G^b=\rho^{-1}(\intZ_2^b)$  and $F\hookrightarrow G^a$,
$F\hookrightarrow G^b$ are the inclusions (see \cite{SW}, p.178).\\

\n{\bf Lemma 3.1:}  \textit{ Let $F$ and $G$ as above. Then}\\

\n{\bf (1)}  \textit{ $F$ has index two in $G^a$ and $G^b$, i.e.
$|G^a/F|=|G^b/F|=2$.}\\

\n{\bf (2)}  \textit{ $G^a$ and $G^b$ are the isotropies of $p_a$
and $p_b$ respectively, i.e. $G^a=G^{\{ p_a\} }$ and
$G^b=G^{\{ p_b\} }$.}
\\

\n{\bf Proof of (1):} Follows from the definitions that  
$G^a/F\cong G^b/F\cong \intZ_2$\\

\n{\bf Proof of (2):} Let $\tilde a\in G$, such that $\rho(\tilde
a)=a$. If $g\in G^a$, then $g=f$ or $g={\tilde a}f$, for some $f\in
F$. Since ${\tilde a},f\in G^{\{p_a\}}$, follows that $g\in G^{\{p_a\}}$.
Conversely, if $g\in G^{\{p_a\}}$,
$gp_a=p_a$, hence $\rho(g)(p_a)=p_a$.
Thus $\rho(g)$ is the identity or $a\in D_{\infty}$. This proves
$G^a=G^{\{p_a\}}$. The proof of $G^b=G^{\{p_b\}}$ is identical. \\

\vspace{0.5cm}

 \n{\bf Proposition 3.2:} \textit{Let $G$ be a
subgroup of a 3-crystallographic group. Suppose that $G$ has a
dihedral action on a line $l$. Then $G$ is isomorphic to  $G^a*_F
G^b$ where the possibilities for $(F, G^a, G^b)$ are  shown in 
table 1.}\\

\begin{table}[h]
\centering
\begin{tabular}{|c||c|c|}
\hline
$F$ & $G^a$ & $G^b$ \\
\hline \hline $ trivial$ & ${\intZ}_{2}$  &
${\intZ}_{2}$ \\
\hline
$\intZ_{2}$ & ${\intZ}_{4}$ & ${\intZ}_{4}$    \\
\cline{2-3}
 & $D_2$ & ${\intZ}_{4}$ \\
\cline{2-3} & $D_2$ & $D_2$ \\
\hline
$\intZ_{3}$ &  ${\intZ}_{6}$ &  ${\intZ}_{6}$     \\
\cline{2-3} & $D_3$ & ${\intZ}_{6}$ \\
\cline{2-3} & $D_3$ & $D_3$ \\
\hline $\intZ_{4}$ & ${\intZ}_{4}\times{\intZ}_{2}$  &
${\intZ}_{4}\times{\intZ}_{2}$ \\
\cline{2-3} & $D_4$ & ${\intZ}_{4}\times{\intZ}_{2}$ \\
\cline{2-3} & $D_4$ & $D_4$ \\
\hline $\intZ_{6}$ & ${\intZ}_{6}\times{\intZ}_{2}$  &
${\intZ}_{6}\times{\intZ}_{2}$ \\
\cline{2-3} & $D_6$ & ${\intZ}_{6}\times{\intZ}_{2}$ \\
\cline{2-3} & $D_6$ & $D_6$ \\
\hline $ D_2$ & $D_2\times{\intZ}_{2}$  &
$D_2\times{\intZ}_{2}$ \\
\cline{2-3} & $D_4$ & $D_2\times{\intZ}_{2}$ \\
\cline{2-3} & $D_4$ & $D_4$ \\
\hline $ D_3$ & $D_3\times{\intZ}_{2}$  &
$D_3\times{\intZ}_{2}$ \\
\hline $ D_4$ & $D_4\times{\intZ}_{2}$  &
$D_4\times{\intZ}_{2}$ \\
\hline $ D_6$ & $D_6\times{\intZ}_{2}$  &
$D_6\times{\intZ}_{2}$ \\
\hline
\end{tabular}
\caption{Dihedral action}
\end{table}

\vspace{0.5cm}

 \n{\bf Remark:} In table 1 we do not give
the inclusions $F\hookrightarrow G^a$ or $F\hookrightarrow G^b$
because it is straightforward to show that if $i_1:F\rightarrow
G^a$, $i_2:F\rightarrow G^a$ are one-to-one, then there is
$\varphi\in Aut(G^a)$ such that $i_2=\varphi\circ i_1$. Analogously for
$G^b$. Hence the isomorphism class of $G^a*_F G^b$ does not
depend on
the particular inclusions.\\

\n{\bf Proof:} Suppose that $G$ acts cocompactly on the line $l$.
By the discussion above, $0\rightarrow F\hookrightarrow
G\stackrel{\rho}\rightarrow \intZ_2^a*\intZ_2^b \rightarrow 0$
and $G=G^a*_{F}G^b$, where\\

\n{\bf (1)} $F$ is trivial or isomorphic to $\intZ_i$, $D_i$,
$i=2,3,4,6$.\\

\n{\bf (2)} $F$ has index two in $G^a$ and $G^b$.\\

\n{\bf (3)} $F=\{g\in G; g|_l=1_l\}$, $G^a=G^{\{p_a\}}$ and
$G^b=G^{\{p_b\}}$ for some $p_a, p_b\in l$, $p_a\neq p_b$.\\

Note that, since $G^a$, $G^b$ fix a point in $\real^3$, $G^a$,
$G^b$ are isomorphic to one of the following groups $\intZ_{i},
\intZ_{i}\times\intZ_{2}, D_{i}, D_{i}\times\intZ_{2}, A_{4},
A_{4}\times\intZ_{2}, S_{4}, S_{4}\times\intZ_{2}, i=2,3,4,6$. By
(1) and (2) above and property (d) (mentioned at the beginning of this
section) we have the following possibilities, shown in
table 2:

\begin{table}[h]
\centering
\begin{tabular}{|c||c|c|}
\hline
$F$ & $G^a$ & $G^b$ \\
\hline \hline $ trivial$ & ${\intZ}_{2}$  &
${\intZ}_{2}$ \\
\hline
$\intZ_{2}$ & ${\intZ}_{4}$ & ${\intZ}_{4}$    \\
\cline{2-3}
 & $D_2$ & ${\intZ}_{4}$ \\
\cline{2-3} & $D_2$ & $D_2$ \\
\hline
$\intZ_{3}$ &  ${\intZ}_{6}$ &  ${\intZ}_{6}$     \\
\cline{2-3} & $D_3$ & ${\intZ}_{6}$ \\
\cline{2-3} & $D_3$ & $D_3$ \\
\hline $\intZ_{4}$ & ${\intZ}_{4}\times{\intZ}_{2}$  &
${\intZ}_{4}\times{\intZ}_{2}$ \\
\cline{2-3} & $D_4$ & ${\intZ}_{4}\times{\intZ}_{2}$ \\
\cline{2-3} & $D_4$ & $D_4$ \\
\hline $\intZ_{6}$ & ${\intZ}_{6}\times{\intZ}_{2}$  &
${\intZ}_{6}\times{\intZ}_{2}$ \\
\cline{2-3} & $D_6$ & ${\intZ}_{6}\times{\intZ}_{2}$ \\
\cline{2-3} & $D_6$ & $D_6$ \\
\hline $ D_2$ & $D_2\times{\intZ}_{2}$  &
$D_2\times{\intZ}_{2}$ \\
\cline{2-3} & $D_4$ & $D_2\times{\intZ}_{2}$ \\
\cline{2-3} & $D_4$ & $D_4$ \\
\cline{2-3} & ${\intZ}_{4}\times{\intZ}_{2}$ & $D_2\times{\intZ}_{2}$ \\
\cline{2-3} & $D_4$ & ${\intZ}_{4}\times{\intZ}_{2}$ \\
\cline{2-3} & ${\intZ}_{4}\times{\intZ}_{2}$ &
${\intZ}_{4}\times{\intZ}_{2}$ \\
\hline $ D_3$ & $D_3\times{\intZ}_{2}$  &
$D_3\times{\intZ}_{2}$ \\
\hline $ D_4$ & $D_4\times{\intZ}_{2}$  &
$D_4\times{\intZ}_{2}$ \\
\hline $ D_6$ & $D_6\times{\intZ}_{2}$  &
$D_6\times{\intZ}_{2}$ \\
\hline
\end{tabular}
\caption{}
\end{table}

\vspace{0.5cm}

 The group ${\intZ}_{4}\times{\intZ}_{2}$ acts by isometries on
$\real^3$ in a unique way (modulo conjugation). That is, the
$\intZ_{4}$ factor fixes  a line $\tilde l$ pointwise and acts on
a plane $\tilde\alpha$ orthogonal to $\tilde l$, by rotations. The
factor $\intZ_{2}$ acts trivially on $\tilde\alpha$ and
 by
reflection on $\tilde l$. Note that the only line in $\real^3$ fixed (as a
set) by $\intZ_4\times\intZ_2$ is $\tilde l$.
If $F=D_2$ and $G^a$ or
$G^b={\intZ}_{4}\times{\intZ}_{2}$,
 then $\intZ_2$ acts trivially on $\alpha$ and by  reflections on the line
$l$. Also $\intZ_4$ acts trivially on $l$  and by rotations
on $\alpha$. But then, by (3) above,
 we should have $D_2\cong \intZ_4$, which is
a contradiction.  Therefore the groups
$({\intZ}_{4}\times{\intZ}_{2})*_{D_2}({\intZ}_{4}\times{\intZ}_{2})$,
$({\intZ}_{4}\times{\intZ}_{2})*_{D_2}D_4$ and
$({\intZ}_{4}\times{\intZ}_{2})*_{D_2}(D_2\times\intZ_2)$ do not
occur geometrically as virtually infinite cyclic subgroups of a
3-crystallographic group $\Gamma$.
This proves the proposition.\\
\vspace{1cm}

\centerline{\large{\bf $H\cong\intZ$}}

\vspace{0,5cm} Then $G$ acts  only by translations  on $l$ and
fits in the exact sequence $$0\rightarrow F\hookrightarrow
G\stackrel{\rho}\rightarrow H\rightarrow 0,$$ with $F$ trivial or
isomorphic to $\intZ_i$ or $D_i$, $i=2,3,4,6$. Thus $G=F*_F\cong
F\rtimes_{\varphi}\intZ$ for some automorphism
$\varphi:F\rightarrow F$
(see \cite{SW}, p.178).\\

To give an action of $\intZ$ on $F$, means to give an automorphism
$\varphi$ of $F$. Since we want to classify the groups
$F*_F\cong F\rtimes_{\varphi}\intZ$ modulo isomorphism,
$\varphi$ must be an outer automorphism, i.e., $\varphi\in$
$Out(F)$. Moreover if $\varphi_1$, $\varphi_2$ are conjugate in
$Out(F)$, $F\rtimes_{\varphi_1}\intZ$ and
$F\rtimes_{\varphi_2}\intZ$ are isomorphic.
A direct calculation shows that:\\

\n{\bf (1)} $Out(F)$ is trivial,\, for $F$ isomorphic to $\intZ_2$ or
$D_3$.\\

\n{\bf (2)} $Out(F)\cong\intZ_2$,\, for $F$ isomorphic to $\intZ_3,
\intZ_4,
\intZ_6,
D_4$ or $D_6.$\\

\n{\bf (3)} $Out(D_2)\cong D_3.$\\

Note that $D_3$ has only three conjugacy classes: $[\, 1\, ], [\,
\varphi\, ], [\, \phi\, ]$, where $\varphi:D_2\rightarrow D_2$ has
order 2 and $\phi:D_2\rightarrow D_2$ has order 3.\\

When does a
group of the form $F\rtimes_{\varphi}\intZ$, with $F$ trivial or
isomorphic to $\intZ_i$ or $D_i$, $i=2,3,4,6$, occur as a virtually
infinite cyclic subgroup of a 3-crystallographic group?\\

To simplify the notation, assume that $l$ is the z-axis. Follows that we can take
$\alpha=l^{\perp}\cong\real^2\cong\real^2\times\{0\}\subset\real^3$.
We can also assume that every $g\in G$ is of the form
$g=(L,ne_3)$, $L\in O(2)\subset O(3)$, $n\in\intZ$, $e_3=(0,0,1)$
(this means $gx=Lx+ne_3$, for all $x\in \real^3$) and the action of
$G$ on $l=\real e_3$ is generated by $\tau_{e_3}$, where
$\tau_{e_3}$ denotes translation by $e_3$. Since $\rho$ is onto,
there is $J\in O(2)\subset O(3)$ such that $(J,e_3)\in G$. Then
the action of $H\cong\intZ$ on $F$ is generated by conjugation by
$g_0=(J,e_3)$, that is, $g_0^{-1}fg_0=J^{-1}fJ=\varphi f$.
 Hence $F\rtimes_{\varphi}\intZ$   occurs
geometrically if there is a $J\in O(2)$ such
that  $J^{-1}fJ=\varphi f$ for all $f\in F$.\\

\n{\bf Proposition 3.3:}  \textit{Let $G$ be a subgroup of a
3-crystallographic group. Suppose that $G$ acts cocompactly by translations on a
line $l$. Then $G$ is isomorphic to $F*_F$, where the possibilities for
$(F,F*_F)$ are shown in table 3:}\\

\begin{table}[h]
\centering
\begin{tabular}{|c|c|}
\hline
$F$ & $F*_F$ \\
\hline \hline
$trivial$ & $\intZ$\\
 \hline
$\intZ_{2}$ & ${\intZ}_{2}\times{{\intZ}}$    \\
\hline
$\intZ_{3}$ & ${\intZ}_{3}\times{{\intZ}}$ \\
\cline{2-2}
 & ${\intZ}_{3}\rtimes_{\varphi}{{\intZ}}$\\
\hline
$\intZ_{4}$ & ${\intZ}_{4}\times{{\intZ}}$\\
\cline{2-2}
 & ${\intZ}_{4}\rtimes_{\varphi}{{\intZ}}$\\
\hline
$\intZ_{6}$ & ${\intZ}_{6}\times{{\intZ}}$\\
\cline{2-2}
  & ${\intZ}_{6}\rtimes_{\varphi}{{\intZ}}$\\
\hline
$D_2$ & $D_2\times{{\intZ}}$\\
\cline{2-2}
  & $D_2\rtimes_{\varphi}{{\intZ}}$\\
\hline
$D_3$ & $D_3\times{{\intZ}}$\\
\hline
$D_4$ & $D_4\times{{\intZ}}$\\
\cline{2-2}
  & $D_4\rtimes_{\varphi}{{\intZ}}$\\
\hline
$D_6$ & $D_6\times{{\intZ}}$\\
\cline{2-2}
  & $D_6\rtimes_{\varphi}{{\intZ}}$\\
\hline
\end{tabular}
\caption{Action by translations}
\end{table}

\vspace{0.5cm}

\n \textit{where $\varphi\neq 1\in$ Out$(F)$, and for $F=D_2$, $\varphi$
has order 2.}\\

\n{\bf Proof:} By the discussion above, we must determine an
automorphism $\varphi:F\rightarrow F$ and  $J\in O(2)$ such that
$J^{-1}fJ=\varphi f$, for all $f\in F$.\\

If $F$ is trivial, then $G=\intZ$\\

If $F=\intZ_i$, $i=3,4,6$, then $Out(\intZ_i)=Aut(\intZ_i)=\{1_F,
\varphi\}$, where $\intZ_{i} = \langle\, t;\, t^{i}=1\, \rangle$
and $\varphi:\intZ_i\rightarrow \intZ_i$, defined by
 $\varphi(t)=t^{i-1}$ and $t\in\intZ^i$ acts on $\real^2$ rotating an angle
 $2\pi /i$. 
Hence, geometrically,  we obtain the groups
$\intZ_i\rtimes_{1_{{\tiny\intZ_i}}}\intZ=\intZ_i\times\intZ$,
$i=2,3,4,6$ by taking $J=1_{\tiny\real^2}$. Also, we obtain
$\intZ_i\rtimes_{\varphi}\intZ$  by taking $
J(x,y)=(x,-y)$, for $i=3,4,6$. \\

Let $F=D_2=\{1, r_x, r_y, -1\}$, where $r_x, r_y$ denote
reflections about the $x$ and $ y$ axes, respectively, and (-1) denotes
rotation by $\pi$. Then
$$Out(D_2)=Aut(D_2)\cong D_3=
\langle \varphi,\, \phi\, ;\,  \phi^3=1,\, \varphi^2=1,\,
(\varphi\phi)^2=1 \rangle.$$

\n Since there are three conjugacy classes  in $D_3$, $[1],
[\varphi]$ and $[\phi]$, there are at most (modulo isomorphism)
three possibilities: $D_2\rtimes_{1_{D_2}}\intZ$,
$D_2\rtimes_{\varphi}\intZ$, and $D_2\rtimes_{\phi}\intZ$. But there is no
$J\in O(2)$ such that $J^{-1}(-1)J\neq (-1)$. Hence
the group $D_2\rtimes_{\phi}\intZ$, does not occur geometrically.
Therefore there are only two possibilities:
$D_2\rtimes_{1_{D_2}}\intZ=D_2\times\intZ$, $D_2\rtimes_{\varphi}\intZ$
which we obtain geometrically taking
$J=1_{\tiny\real^2}$ and
$J(x,y)=(y,x)$, respectively.\\

If $F$ is isomorphic to $\intZ_2$ or $D_3$, then $Out(F)=1$ and we
obtain the  groups $\intZ_2\times\intZ$ and
$D_3\times\intZ$.\\

If $F$ is isomorphic to $D_4$ or $D_6$, then $F= \langle s,\, t; t^i=1,\,
s^2=1,\, (st)^2=1
\rangle $,\, $i=4,6$,\, and we have $Out(D_4)\cong Out(D_6)\cong\{1,
\varphi\}$, where $\varphi(s)=st$ and $\varphi(t)=t$,\,
$i=4,6$. Hence there are two possibilities:
$D_i\rtimes_{1_{D_i}}\intZ=\intZ_i\times\intZ$ and
$D_i\rtimes_{\varphi}\intZ$, $i=4,6$, and we can show that they
occur geometrically by taking  $J=1_{\tiny\real^2}$ or 
equal to a rotation by an angle $\pi /i$, for $i=4,6$. This proves the proposition.

\vspace{1cm}

\section{ Isotropy in  dimension three}

\hspace{0.5cm} In this section we calculate the possible isotropy groups
$\Gamma^\sigma=\{\gamma\in \Gamma;\gamma\sigma=\sigma\}$, for
an open cell $\sigma$ in the   universal
 $(\Gamma,VC(\Gamma))$-space
$A$.

We will need the following lemma.\\

\n{\bf Lemma 4.1} \textit{Let $G$ be a virtually infinite cyclic
subgroup of a 3-crystallographic group $\Gamma$. If $G$ leaves
invariant two different lines, then $G$ is trivial or isomorphic to
$\intZ,
D_{\infty}, \intZ\times\intZ_2, D_{\infty}\times\intZ_2$}.\\

\n{\bf Proof:} Suppose that $Gl=l$ and $G\tilde l =\tilde l$,
with $l\neq\tilde l$, and $G$ acts cocompactly on $l$ (see prop. 1.1). Hence
$l\parallel\tilde l$ (if $l,\tilde l$ are not contained in a
plane, $G$ is trivial). Thus $G$ leaves invariant the plane
$\alpha$ that contains $l, \tilde l$. Moreover $G$ leaves
invariant all lines in $\alpha$ which are parallel to $l$. Note
that if  $g\in G$ acts trivially on $l$, then $g$ acts trivially
on $\alpha$. Let $\beta$ be the plane that contains $l$ and is
orthogonal to the plane $\alpha$. Then $G$ leaves invariant
$\beta$. Note that if $g\in G$ acts trivially on $\beta$, then $g$
acts trivially on $\real^3$, i.e., $g$ is the identity. Hence $G$
is a virtually infinite cyclic group acting faithfully on $\real^2$.
Therefore, $G$ is isomorphic to $ \intZ, D_{\infty},
\intZ\times\intZ_2, D_{\infty}\times\intZ_2$ (see \cite{P}).
This proves the lemma.\\

Recall that every (open) cell in $A - Z=(C^3\times\real^3)- Y$ is a 
product of (open) simplices $\sigma_C \times\sigma_R$ (see sect. 2).

\vspace{1cm}

\centerline {\large{\bf Isotropy of 0-cells}}

\vspace{1cm}

 Let $\sigma^{0}$ be a 0-cell in $A$ and
 $\gamma\in\Gamma$.
We have two cases:\\

\n{\bf First case:}  $\sigma^{0}\subset A-Z$.\\

Recall that $A-Z$ is $\Gamma$-homeomorphic to
$(C^3\times\real^3)-Y$. Then
$\sigma^{0}=\sigma^{0}_C\times\sigma^{0}_R$.
 If $\gamma\sigma^0=\sigma^0$,
 we have
that $\gamma\sigma^0_C=\sigma^0_C$ and
$\gamma\sigma^0_R=\sigma^0_R$. Since $\gamma$ fixes the point
$\sigma^0_R$ in $\real^3$, $\gamma$ is a rotation. Thus
$\Gamma^{\sigma}\, <\, O(3)$. Hence we have
two possibilities:\\

 \n{\bf (1)} If $\sigma^0_C$ is not  the vertex $[ 0 ]$ of
$C^3$ (the vertex does not determine any direction), $\gamma$
fixes a direction in $\real^3$, the direction determined by
$\sigma^{0}_C\neq[ 0 ]$. Therefore $\gamma$ leaves invariant a
line $l$ in $\real^3$ (the line that passes through
$\sigma^{0}_R$ and has direction determined by $\sigma^{0}_C$). By
(d) at the beginning of section 3,  follows that $\Gamma^{\sigma^0}$ is trivial or
isomorphic to one of groups $D_i, \intZ_i, D_i\times \intZ_2,
\intZ_i\times \intZ_2, \,
i=2,3,4,6$.\\

\n{\bf (2)} If $\sigma^0_C$ is the vertex [0] of $C^3$,
$\sigma^0_C$ does not determine any direction in $\real^3$. In
this case the only thing we can say is that $\gamma$  fixes a
point in $\real^3$. Therefore, by (b) at the beginning of section 3, $\Gamma^{\sigma^0}$ is
trivial or isomorphic to one of groups $D_i, \intZ_i, A_4, S_4,
D_i\times \intZ_2, \intZ_i\times \intZ_2, A_4\times \intZ_2,
S_4\times \intZ_2,\,
i=2,3,4,6$.\\

\n{\bf Second case:}  $\sigma^{0}\subset Z$.\\

Remark that $Z=\bigcup_{h\in H}Z_{h}$. Then $\sigma^{0}\subset
Z_{h}$ for some $h\in H$.
 If $\gamma\sigma^{0}=\sigma^{0}$, $\gamma$ fixes
a  point in $Z_{h}$. Recall  that a point in $Z_{h}$ is obtained
by
 collapsing  a line in $C^3\times\real^3$. Let $l$ be the line the collapses
 to $\sigma^{0}$. Then $\Gamma^{\sigma^{0}}=\Gamma^{(l)}$.\\

\vspace{1cm}

\centerline{\large{\bf Isotropy of 1-cells}}

\vspace{1cm}

Let $\sigma^{1}$ be an open 1-cell in $A$, and $\gamma\in\Gamma$.
We have two cases:\\

\n{\bf First case:}  $\sigma^{1}\subset A-Z$.\\

We have two possibilities:
$\sigma^{1}=\sigma^{0}_{C}\times\sigma^{1}_{R}$ or
$\sigma^{1}=\sigma^{1}_{C}\times\sigma^{0}_{R}$.\\

\n{\bf (1)} Suppose
$\sigma^{1}=\sigma^{0}_{C}\times\sigma^{1}_{R}$. If
$\gamma\sigma^{1}=\sigma^{1}$ then
$\gamma\sigma_{R}^{1}=\sigma_{R}^{1}$. Hence $\gamma$ fixes   a
line pointwise. Therefore $\gamma$ leaves invariant any plane
orthogonal to this line and acts by rotation or reflection on these planes. By
the crystallographic restriction, we have that
$\Gamma^{\sigma^{1}}$ is trivial or isomorphic to one of the
groups $D_{i}, \intZ_{i},
i=2,3,4,6.$\\

\n{\bf (2)} Suppose
$\sigma^{1}=\sigma^{1}_{C}\times\sigma^{0}_{R}$. Thus if
$\gamma\sigma^{1}=\sigma^{1}$ then
$\gamma\sigma_{R}^{0}=\sigma_{R}^{0}$ and
$\gamma\sigma_{C}^{1}=\sigma_{C}^{1}$. Hence  $\gamma$ fixes one
point in $\real^3$ and at least, one direction in $\real^{3}$ (a
direction determined by a point in $\sigma_{C}^{1}$, different
from the vertex $[0]$). By (d) of section 3, $\Gamma^{\sigma^{1}}$ is
trivial or isomorphic to one of the groups $ D_i, \intZ_i,
D_{i}\times\intZ_{2}, \intZ_{i}\times\intZ_{2},
i=2,3,4,6.$\\

\n{\bf Second case} $\sigma^{1}\subset Z$.\\

Since
 $Z=\bigcup_{h\in H}Z_{h}$, $\sigma^{1}\subset
Z_{h}$ for some $h\in H$.
 Then if $\gamma\sigma^{1}=\sigma^{1}$, $\gamma$ fixes
pointwise a line  in $Z_{h}$. Hence, $\Gamma^{\sigma^{1}}$ fixes,
at least, two points in $Z_h$. Therefore $\Gamma^{\sigma^{1}}$
leaves invariant, at least, two lines in $\real^3$. Then Lemma 4.1
implies that $\Gamma^{\sigma^{1}}$ is trivial or isomorphic to $
\intZ,
D_{\infty}, \intZ\times\intZ_2, D_{\infty}\times\intZ_2$.\\

\vspace{1cm}

\centerline {\large{\bf Isotropy of 2-cells}}

\vspace{1cm}

Let $\sigma^{2}$ be an open 2-cell in $A$ and $\gamma\in
\Gamma^{\sigma^{2}}$.
We have two cases:\\

\n{\bf First case:}  $\sigma^{2}\subset A-Z$.\\

 We have three
possibilities: $\sigma^{2}=\sigma^{0}_{C}\times\sigma^{2}_{R}$,
$\sigma^{2}=\sigma^{2}_{C}\times\sigma^{0}_{R}$ or
$\sigma^{2}=\sigma^{1}_{C}\times\sigma^{1}_{R}.$
\\

\n{\bf (1)} Suppose
$\sigma^{2}=\sigma^{0}_{C}\times\sigma^{2}_{R}$. If
$\gamma\sigma^{2}=\sigma^{2}$ then
$\gamma\sigma_{R}^{2}=\sigma^{2}_{R}$, and follows that
 $\gamma$ fixes a plane in $\real^3$ pointwise.
Hence, $\Gamma^{\sigma^{2}}$ is trivial or isomorphic to $\intZ_{2}$.\\

\n{\bf (2)} Suppose
$\sigma^{2}=\sigma^{2}_{C}\times\sigma^{0}_{R}$. Thus if
$\gamma\sigma^{2}=\sigma^{2}$ then $
\gamma\sigma_{C}^{2}=\sigma^{2}_{C}$ and
$\gamma\sigma_{R}^{0}=\sigma^{0}_{R}$. Hence $\gamma$ fixes
infinitely many directions in $\real^3$, and fixes a point in
$\real^3$. Then $\gamma$ leaves invariant at least a plane in
$\real^3$ and acts on it (at most) by rotation by $\pi$. Also,
$\gamma$ acts (at most) by reflection on the line orthogonal  to
this plane. Therefore, $\Gamma^{\sigma^{2}}$ is trivial or
isomorphic to $\intZ_{2}$ or
$\intZ_{2}\times\intZ_{2}$.\\

 \n{\bf (3)}
Suppose $\sigma^{2}=\sigma^{1}_{C}\times\sigma^{1}_{R}$. Thus if
$\gamma\sigma^{2}=\sigma^{2}$ then
$\gamma\sigma^{1}_{R}=\sigma^{1}_{R}$. Hence $\gamma$ fixes a line
in $\real^3 $ pointwise. Then $\gamma$ acts on the
plane orthogonal  to this line (if $\sigma^{1}_{C}$ determines the
same direction as $\sigma^{1}_{R}$). Therefore,
$\Gamma^{\sigma^{2}}$ is trivial or isomorphic to one of $D_{i}$
or $\intZ_{i}$, $i=2,3,4,6$.\\

Note that if $\Gamma^{\sigma^{2}}$ is isomorphic to $D_{i}$ or
$\intZ_{i}$, $i=2,3,4,6$, then $\sigma^{1}_{C}$ is a ray cell, and
the direction that $\sigma^{1}_{C}$ determines is exactly the
direction of $\sigma^{1}_{R}$. In this case, we say that
$\sigma^{2}$ is a {\it special 2-cell}.

\vspace{0.5cm}

\n{\bf Second case} $\sigma^{2}\subset Z$.\\

Since  $Z=\bigcup_{h\in H}Z_{h}$, then $\sigma^2\subset Z_h$ for
some $h\in H$. If $\gamma\sigma^{2}=\sigma^{2}$, $\gamma$ fixes
$Z_{h}$ pointwise. Hence $\gamma$ leaves invariant all lines
parallel to $h$. Since $\Gamma^{\sigma^{2}}$ acts cocompactly on
$l_h$, follows that $\Gamma^{\sigma^{2}}$  acts faithfully on
$l_h$. Therefore, $\Gamma^{\sigma^{2}}$ is trivial or
isomorphic to $D_{\infty}$ or  $\intZ$.\\

\vspace{1cm}

\centerline {\large{\bf Isotropy of 3-cells}}

\vspace{1cm}

 Let $\sigma^{3}$ be an open 3-cell in $A$ and
 $\gamma\in\Gamma^{\sigma^{3}}$.
Since $Z=\bigcup_{h\in H}Z_{h}$ and $Z_h\cong\real^2$, we have that
$\sigma^{3}\subset A-Z$. Hence, there are four possibilities:
 $\sigma^{3}=\sigma^{0}_{C}\times\sigma^{3}_{R}$,
$\sigma^{3}=\sigma^{1}_{C}\times\sigma^{2}_{R}$,
$\sigma^{3}=\sigma^{2}_{C}\times\sigma^{1}_{R}$ or
$\sigma^{3}=\sigma^{3}_{C}\times\sigma^{0}_{R}$.\\

\n{\bf (1)} Suppose $\sigma^{3}=\sigma^{0}_{C}\times\sigma^{3}_{R}$. If
$\gamma\sigma^{3}=\sigma^{3}$, then $\gamma\sigma^{3}_R=\sigma^{3}_R$.
 Hence $\gamma$ fixes the whole $\real^3$ pointwise.
Therefore,
$\Gamma^{\sigma^{3}}$ is trivial.\\

\n{\bf (2)} Suppose $\sigma^{3}=\sigma^{1}_{C}\times\sigma^{2}_{R}$. If
$\gamma\sigma^{3}=\sigma^{3}$, then $\gamma\sigma^{2}_R=\sigma^{2}_R$
and $\gamma\sigma^{1}_C=\sigma^{1}_C$. Hence $\gamma$ fixes pointwise a
plane $\alpha$ and fixes a direction in $\real^3$. Then $\gamma$ acts
trivially
or by reflection on the line orthogonal to $\alpha$. Therefore,
$\Gamma^{\sigma^{3}}$ is trivial or isomorphic to $\intZ_2$.\\

\n{\bf (3)} Suppose  $\sigma^{3}=\sigma^{2}_{C}\times\sigma^{1}_{R}$. If
$\gamma\sigma^{3}=\sigma^{3}$, then $\gamma\sigma^{1}_R=\sigma^{1}_R$
and $\gamma\sigma^{2}_C=\sigma^{2}_C$. Hence $\gamma$ fixes a line
pointwise
 and fixes infinitely many directions in $\real^3$. Then $\gamma$ leaves
invariant the plane orthogonal to the line
fixed by $\gamma$ and acts on it (at most) by
rotation by $\pi$.  Therefore,
$\Gamma^{\sigma^{3}}$ is trivial or isomorphic to $\intZ_2$.\\

\n{\bf (4)} Suppose
$\sigma^{3}=\sigma^{3}_{C}\times\sigma^{0}_{R}$. If
$\gamma\sigma^{3}=\sigma^{3}$, then
$\gamma\sigma^{0}_R=\sigma^{0}_R$ and
$\gamma\sigma^{3}_C=\sigma^{3}_C$. Hence $\gamma$ fixes a point
and fixes all directions in $\real^3$. Therefore,
$\Gamma^{\sigma^{3}}$ is trivial or isomorphic to $\intZ_2$.\\

\vspace{1cm}

\centerline {\large{\bf Isotropy of 4-cells}}

\vspace{1cm}

 Let $\sigma^{4}$ be an open 4-cell in $A$ and
 $\gamma\in\Gamma$.
Since $Z=\bigcup_{h\in H}Z_{h}$ and $Z_h\cong\real^2$, we have
$\sigma^{4}\subset A-Z$. Hence, there are three possibilities:
 $\sigma^{4}=\sigma^{1}_{C}\times\sigma^{3}_{R}$,
$\sigma^{4}=\sigma^{2}_{C}\times\sigma^{2}_{R}$ or
$\sigma^{4}=\sigma^{3}_{C}\times\sigma^{1}_{R}$.\\

\n{\bf (1)} Suppose $\sigma^{4}=\sigma^{1}_{C}\times\sigma^{3}_{R}$. If
$\gamma\sigma^{4}=\sigma^{4}$, then $\gamma\sigma^{3}_R=\sigma^{3}_R$.
 Hence $\gamma$ fixes whole $\real^3$ pointwise. Therefore,
$\Gamma^{\sigma^{4}}$ is trivial.\\

\n{\bf (2)} Suppose
$\sigma^{4}=\sigma^{2}_{C}\times\sigma^{2}_{R}$. If
$\gamma\sigma^{4}=\sigma^{4}$, then
$\gamma\sigma^{2}_R=\sigma^{2}_R$ and
$\gamma\sigma^{2}_C=\sigma^{2}_C$. Hence $\gamma$ fixes  a plane
pointwise. Therefore, $\Gamma^{\sigma^{4}}$ is trivial or
isomorphic to
$\intZ_2$.\\

\n{\bf (3)} Suppose  $\sigma^{4}=\sigma^{3}_{C}\times\sigma^{1}_{R}$. If
$\gamma\sigma^{4}=\sigma^{4}$, then $\gamma\sigma^{1}_R=\sigma^{1}_R$
and $\gamma\sigma^{3}_C=\sigma^{3}_C$. Hence $\gamma$ fixes  a
line pointwise and
all directions in $\real^3$.
  Therefore,
$\Gamma^{\sigma^{4}}$ is trivial.\\

\vspace{1cm}

\centerline {\large{\bf Isotropy of 5-cells}}

\vspace{1cm}

 Let $\sigma^{5}$ be an open 5-cell in $A$ and
 $\gamma\in\Gamma$.
Since $Z=\bigcup_{h\in H}Z_{h}$ and $Z_h\cong\real^2$, we have
$\sigma^{5}\subset A-Z$. Hence, there are two possibilities:
 $\sigma^{5}=\sigma^{2}_{C}\times\sigma^{3}_{R}$ or
$\sigma^{5}=\sigma^{3}_{C}\times\sigma^{2}_{R}$.\\

\n{\bf (1)} Suppose $\sigma^{5}=\sigma^{2}_{C}\times\sigma^{3}_{R}$. If
$\gamma\sigma^{5}=\sigma^{5}$, then $\gamma\sigma^{3}_R=\sigma^{3}_R$.
 Hence $\gamma$ fixes  whole $\real^3$ pointwise.
Therefore,
$\Gamma^{\sigma^{3}}$ is trivial.\\

\n{\bf (2)} Suppose  $\sigma^{5}=\sigma^{3}_{C}\times\sigma^{2}_{R}$. If
$\gamma\sigma^{5}=\sigma^{5}$, then $\gamma\sigma^{2}_R=\sigma^{2}_R$
and $\gamma\sigma^{3}_C=\sigma^{3}_C$. Hence $\gamma$ fixes  a
plane pointwise and
all directions in $\real^3$.
  Therefore,
$\Gamma^{\sigma^{4}}$ is trivial.\\

\vspace{1cm}

\centerline {\large{\bf Isotropy of 6-cells}}

\vspace{1cm}

 Let $\sigma^{6}$ be an open 6-cell in $A$ and
 $\gamma\in\Gamma$.
Since $Z=\bigcup_{h\in H}Z_{h}$ and $Z_h\cong\real^2$, we have
$\sigma^{6}\subset A-Z$. Hence
$\sigma^{6}=\sigma^{3}_{C}\times\sigma^{3}_{R}$. If
$\gamma\sigma^{6}=\sigma^{6}$, then $\gamma\sigma^{3}_R=\sigma^{3}_R$.
Hence, $\Gamma^{\sigma^{6}}$ is trivial.\\

\vspace{1cm}

\section{ Calculation of $H_{i}(A/\Gamma,{\cal P}_{*}(\rho))$}

\hspace{0.5cm} Recall that a cell $\tilde\sigma$ in $A/\Gamma$
corresponds to an orbit
$\Gamma\sigma, \, \sigma\in A.$
In what follows we will use the same notation for a cell $\sigma$
and its orbit $\Gamma\sigma$.
We know that there is a spectral sequence with
\vspace{0.1cm}
${\cal
E}^{2}_{p,q}=H_{p}(A/\Gamma,\pi_{q}({\cal P}_{*}(\rho)))$ (see Prop. 1.2 and 1.4)
which abuts to $H_{p+q}(A/\Gamma,{\cal P}_{*}(\rho))$. We are interested in
the case  $p+q=-1$. By Prop. 1.3, \cite{Ca} and \cite{FJ2}, $\pi_{q}({\cal
P}_{*}(\rho))=0$, if $q+2\leq-2$. Hence
$H_{p}(A/\Gamma,\pi_{q}({\cal P}_{*}(\rho)))=0$ if $q+2\leq-2$.
Then, the possible non zero terms of the spectral sequence with
$p+q=-1$ are ${\cal E}^{2}_{0,-1}=H_{0}(A/\Gamma,\pi_{-1}({\cal
P}_{*}(\rho)))$, ${\cal
E}^{2}_{1,-2}=H_{1}(A/\Gamma,\pi_{-2}({\cal P}_{*}(\rho)))$ and
${\cal E}^{2}_{2,-3}=H_{2}(A/\Gamma,\pi_{-3}({\cal
P}_{*}(\rho)))$.\\

\vspace{1cm}

\centerline{\large{\bf Calculation of the term ${\cal E}^{2}_{0,-1}$.}}

\vspace{1cm}

Let $\sigma_{i}^{0}$ and $\sigma^{1}_{j}$ denote the 0-cells and 1-cells of $A$.
Consider the associated cellular chain complex $...\leftarrow
C_{0}\stackrel{\partial}\longleftarrow C_{1}
\leftarrow ...$,
 for ${\cal E}^{2}_{0,-1}$, where $C_{0}=
\bigoplus_{i}Wh(\Gamma^{\sigma_{i}^{0}})$, $C_{1}=
\bigoplus_{j}Wh(\Gamma^{\sigma^{1}_{j}})$.
By
the calculations of Section 4 (isotropy of 1-cells),
$\Gamma^{\sigma^{1}}$ is trivial or isomorphic to one of groups
$\intZ\times\intZ_{2}, D_{\infty}\times\intZ_{2}, D_{i}, \intZ_i,
\intZ, D_{\infty},
 \intZ_i\times\intZ_2, D_i\times\intZ_2, i=2,3,4,6.$
In \cite{W}, Whitehead proves that $Wh(F)=0$ if $F$ is trivial or isomorphic
to
$\intZ_i, \, i=2, 3, 4$.
 In \cite{Ca}, Carter proves that $Wh(\intZ_6)=0$. In \cite{B2},
Bass proves that $Wh(\intZ)=Wh(D_{\infty})=0$. In \cite{P}, Pearson proves
that
$Wh(\intZ\times\intZ_2)=Wh(D_{\infty}\times\intZ_2)=0$. In \cite{B1}, Bass
proves
that $Wh(F)=\intZ^{r-q}+SK_1(\intZ[F])$, when $F$ is a finite group, where
$r$ is the number of  real representations irreducible of $F$ and
$q$ is the  number of  rational representations irreducible. By 
Theorems 14.1 and 14.2 of \cite{O}, we have $SK_1(\intZ[F])=0$, if $F$ is isomorphic to
$D_i, D_i\times\intZ_2, \intZ_4\times\intZ_2, \intZ_6\times\intZ_2$,
$i=2,3,4,6$. A direct calculation shows that $r=q$ for these groups.
Then $Wh(\Gamma^{\sigma^{1}})=0$. Hence $C_{1}=
\bigoplus_{j}Wh(\Gamma^{\sigma^{1}_{j}})=0$.
Therefore,
$${\cal
E}^{2}_{0,-1}=H_{0}(A/\Gamma,\pi_{-1}({\cal P}_{*}(\rho)))=
\bigoplus_{i}Wh(\Gamma^{\sigma_{i}^{0}}).$$

The  following lemmata show that the term ${\cal E}^{2}_{0,-1}$
can be written in terms of the maximal virtually cyclic
subgroups of $\Gamma$.\\

 \n{\bf Lemma 5.1:}
\textit{ Let $l\in {\cal L}$ and suppose that $\Gamma^{(l)}$ is a
 virtually cyclic infinite subgroup of $\Gamma$.
If ${\cal L}^{\Gamma^{(l)}}=\{l\}$,
then $\Gamma^{(l)}$ is a  maximal virtually finite cyclic
subgroup of $\Gamma$.}\\

\n{\bf Proof:} Suppose that $G$ is virtually infinite cyclic  and
that $\Gamma^{(l)}<G$. Since $G$ leaves invariant at least one
line $\tilde l$, we have $\Gamma^{(l)}<G<\Gamma^{(\tilde l)}$.
Hence $ l=\tilde l$. This proves the lemma.\\

 \n{\bf Lemma 5.2:} \textit{ Let $p\in \real^3$. If
$\Gamma^{p}$ fixes only $p$, then
$\Gamma^{p}$ is a maximal finite subgroup of $\Gamma$.}\\

\n{\bf Proof:} Suppose that $F$ is a finite subgroup of $\Gamma$ and that
$\Gamma^{p}<F$. Since $F$ is finite, it fixes at least a point
${\tilde p}$ in $\real^3$, then
$\Gamma^{p}<F<\Gamma^{\tilde p}$. Hence
${p}={\tilde p}$. This proves the lemma. \\

 \n{\bf Lemma 5.3:} \textit{Let $G$ be a
virtually infinite cyclic subgroup of $\Gamma$. If $Wh(G)\neq 0$,
then ${\cal L}^{G}=\{l\}$} or, equivalently,
$G=\Gamma^{l}$, for some unique line $l$.\\

\n{\bf Proof:} If ${\cal L}^{G}\supset\{l,\tilde l\}$,
$l\neq\tilde l$, by Lemma 4.1, $G$ is  isomorphic to $ \intZ,
D_{\infty}, \intZ\times\intZ_2, D_{\infty}\times\intZ_2$. Follows
that $Wh(G)=0$ (see \cite{P}). This proves the
lemma.\\

 Recall  that $\Gamma^z$, $z\in Z$
is a  virtually infinite cyclic group (because a point in $Z$ is
obtained by collapsing a line in $C^3\times\real^3$, see prop. 1.1). Hence,
lemmata 5.1 and 5.3 imply that if
$Wh(\Gamma^z)\neq 0$, then $\Gamma^z$ is maximal virtually infinite
cyclic. Therefore, we can write
$${\cal E}^{2}_{0,-1}= \left[\bigoplus_{H\in
F(\Gamma)}Wh(H)\right]\oplus\left[\bigoplus_{H\in
{VC}_{\infty}(\Gamma)}Wh(H)\right],$$
where  $F(\Gamma)$ and ${VC}_{\infty}(\Gamma)$ denote the sets of conjugacy
classes of maximal finite subgroups of $\Gamma$ and maximal virtually infinite 
cyclic subgroups of $\Gamma$, respectively.\\

As mentioned above, the Whitehead groups of $\intZ_i$, $D_i$, 
$\intZ_i \times\intZ_2$, $D_i \times\intZ_2$, $i=2,3,4,6$
vanish. We show now that the Whitehead groups of the other finite
subgroups of 3-crystallographic groups also vanish.
The groups $A_4, S_4, A_4\times\intZ_2, S_4\times\intZ_2$, act on
$\real^3$ fixing only a point (see (d) at the beginning of Section 3). 
Then, by Lemma 5.2, these groups are maximal. By \cite{O}, Theorem 14.1, 
follows that $Wh(S_4)=Wh( S_4\times\intZ_2)=0$. By \cite{O}, Theorem 14.6, 
follows that $SK_1(\intZ [A_4])=0$. By \cite{B1}, $Wh(A_4)=\intZ^{r-s}\oplus
SK_1(\intZ [A_4])$. A direct calculation shows that $r=q$. Hence
$Wh(A_4)=0$.\\

The proof of the following lemma was suggested to us by F.T. Farrell.\\

\n{\bf Lemma  5.4:} \textit{ If $F=A_{4}\times\intZ_{2}$, then
$Wh(F)=0$.}\\

\n{\bf Proof:} Recall that $Wh(F)=\intZ^y\oplus SK_1(\intZ [F])$.
Since the elementary subgroups of $A_{4}\times\intZ_{2}$ are either
cyclic or the sum of cyclic subgroups of order 2, by \cite{O}, examples 1 and
2 on page 14, Theorem 5.3 and page 7, we have that $SK_1(\intZ [F])=0$. 
Now the rank of the torsion-free part
of $Wh(F)$ is $y=r-q$, where $r$ is the number of conjugacy classes of
sets $\{ g, g^{-1}\}$ and $q$ is the number of
conjugacy classes of cyclic subgroups of $F$ (see \cite{M}).
We have six cyclic subgroups  (modulo conjugation) of
$A_{4}\times\intZ_{2}$: the  trivial subgroup, three subgroups
isomorphic to  $\intZ_2$, one subgroup isomorphic to $\intZ_{3}$
and one subgroup isomorphic to $\intZ_{6}$. On the other hand,  
in $A_4$ there are exactly three elements of order 2 (all of them are conjugate)
and eight elements of order 3.
Let $a,b\in A_4$ with $a^3=b^3=1$ and $a,b$ non  trivial. A direct
calculation shows that either $a=b^2$ or $a$ and $b$ are
conjugate. Hence there are six  of conjugacy classes of sets $\{
g, g^{-1}\}$ in $A_{4}\times\intZ_{2}$.
Therefore,
 $r=q$ and
$Wh(A_{4}\times\intZ_{2})=0$. This proves the lemma.\\

Follows that $Wh(F)=0$, if $F$ is a finite subgroup of a
3-crystallographic group. Hence we
can write:
$$Wh(\Gamma)=\bigoplus_{G\in
{VC}_{\infty}(\Gamma)}Wh(G).$$\\

\hspace{1cm}

\centerline{\large{\bf Calculation of the term ${\cal E}^{2}_{2,-3}$.}}

\vspace{0,5cm}

Consider the associated cellular chain complex $...\leftarrow
C_{1}\stackrel{\partial}\longleftarrow C_{2}\leftarrow...$, for
${\cal E}^{2}_{2,-3}$, where
$C_{1}=\bigoplus_{j}K_{-1}(\intZ[\Gamma^{\sigma_{j}^{1}}])$ and
$C_{2}=\bigoplus_{i}K_{-1}(\intZ[\Gamma^{\sigma_{i}^{2}}])$.
By Section 4 (isotropy of 2-cells)
$\Gamma^{\sigma_{i}^{2}}$  is trivial or isomorphic to  $\intZ_i, D_i,
\intZ$
or $D_{\infty}$. By \cite{B2}, we have that $K_{-1}(\intZ [G])=0$, if $G$
is isomorphic to
  $D_2, \intZ$, $D_{\infty}$ or $\intZ_i$,
$i=2,3,4$, and that
$K_{-1}(\intZ[\intZ_6])=\intZ$. By \cite{P}, we have that $K_{-1}(\intZ
[G])=0$
if $G$ is isomorphic to $D_3$ or $D_4$ and that $K_{-1}(\intZ
[D_6])=\intZ$.
Hence, it is enough to study the case in which $\sigma^{2}$ is a special 2-cell 
with isotropy $D_6$ (in the other cases $K_{-1}$ vanishes), i.e.,
$\sigma^{2}=\sigma_{C}^{1}\times\sigma^{1}_{R}$, where
$\sigma_{C}^{1}$ is a ray cell that determines the same direction as 
$\sigma^{1}_{R}$ and the isotropy of $\sigma^{1}_{R}$ is $D_6$.\\

In what follows the bar denotes ``closure''.\\

\n{\bf Lemma 5.5:} \textit{Let
$\sigma^{2}=\sigma_{C}^{1}\times\sigma^{1}_{R}$,
${\tilde\sigma}^{2}=\tilde\sigma_{C}^{1}\times\tilde\sigma^{1}_{R}$
be special 2-cells in $A-Z$. If
\vspace{0,1cm}
$\overline{\sigma^{2}}\cap{\overline{\tilde\sigma^{2}}}$  is an
1-cell of the form
$\overline{\sigma_{C}^{0}\times\sigma^{1}_{R}}$, then
$\sigma^{1}_{R}=\tilde\sigma^{1}_{R}$, and
$\sigma_{C}^{1}$, $\tilde\sigma_{C}^{1}$ lie in the same ray in
$C^3$.}\\

\n{\bf Proof:}  Since
$\overline{\sigma^{2}}\cap\overline{{\tilde\sigma}^{2}}
=\overline{\sigma_{C}^{0}\times\sigma^{1}_{R}}$
we have $\sigma_{R}^{1}=\tilde\sigma^{1}_{R}$. Hence
$\sigma_{C}^{1}$ and $\tilde\sigma_{C}^{1}$ determine the same
direction as  $\sigma^{1}_{R}=\tilde\sigma^{1}_{R}$. Therefore
$\sigma_{C}^{1}$ and $\tilde\sigma_{C}^{1}$ lie in the
same ray. This proves the lemma.\\

It follows from the lemma above that the only special 2-cell that
\vspace{0,1cm}
contains
a one cell of the form $[0]\times\sigma^{1}_{R}$ is the
special 2-cell $\overline{\sigma_{C}^{1}\times\sigma^{1}_{R}}$, where
$\sigma_{C}^{1}$ is the ray cell that contains the cone point $[0]$
and determines the same direction as $\sigma_{R}^{1}$.\\

We will prove that $\partial:C_{2}\rightarrow C_{1}$ is injective.
 It will then
follow that ${\cal E}^{2}_{2,-3}=0$.\\

Fix $\sigma^{1}_{R}$  and enumerate all ray cells
$\sigma^{1}_{C}$  which determine the same direction as
$\sigma_{R}^{1}$. Then we have the family of 2-cells
$\{\sigma^{2}_{i}=(\sigma_{C}^{1})_{i}\times\sigma^{1}_{R}\}$.
Note that $\cup_{i}(\overline{\sigma^{1}_{C}})_i$ is the ray $r\cong
[0,\infty)$, determined by the direction of $\sigma^{1}_{R}$.
Denote by $l$  the line that contains $\sigma_R^1$. After
reindexing the cells, we can identify each $(\sigma_{C}^{1})_{i}$
with the interval $(i,i+1)$. Therefore $\sigma^{2}_{i}=
(i,i+1)\times\sigma_{R}^{1}$.\\

From the definition of the relation $\cong$ over $C^{3}\times\real^{3}$
(Section 2) we get that we have to collapse at most one cell of
the form $\{n\}\times\sigma^{1}_{R}$ to one point, where $n>0$
depends on $\sigma^{1}_{R}$. In fact, if $l$ is parallel to
$l_h\in{\cal L}_h\in\Lambda_n$,
for some $h\in H$, then $c(l)=[\frac{nh}{||h||}]$.
Hence $\{n\}\times\sigma^{1}_{R}$ collapses to one point. Therefore $n$
depends on $\sigma^{1}_{R}$. We write $n=n(\sigma^{1}_{R})$. Then,
either $\cup_i\overline{\sigma^{2}_{i}}=[0,\infty)\times\overline{\sigma^{1}_{R}}$, or
$\cup_i\overline{\sigma^{2}_{i}}$ is obtained from
$[0,\infty)\times\overline{\sigma^{1}_{R}}$ collapsing
$\{n\}\times\overline{\sigma^{1}_{R}}$ to a point $\{n\}$ ( if $l\parallel
l_h$, for some $h$), where $n=n(\sigma^{1}_{R})>0$.
\\

\n{\bf Lemma 5.6:} \textit{$\partial:C_2\longrightarrow C_1$ is
injective.}\\

\n{\bf Proof:} Let $\sum_{k}m_{k}\sigma^{2}_{k}$ $\in C_2$ with $\partial 
(\sum_{k}m_{k}\sigma^{2}_{k})=0$.
Since $m_{k}\in {K}_{-1}(\intZ[\Gamma^{\sigma_{k}^{2}}])$, we
can suppose that $\sigma_{k}^{2}$ is a special 2-cell,
that is, it is of the form $\sigma_{C}^{1}\times\sigma^{1}_{R}$,
with $\sigma_{C}^{1}$ a ray cell that determines the same
direction as $\sigma^{1}_{R}$ and the isotropy of $\sigma^{1}_{R}$ is
$D_6$ (in the other cases ${K}_{-1}(\intZ[\Gamma^{\sigma_{k}^{2}}])=0$). 
Thus we can write

$$\sum_{k}m_{k}\sigma^{2}_{k}=\sum_{\sigma^{1}_{R}}\left(
\sum_{i}n_i((\sigma^1_C)_i\times\sigma_R^1)\right).$$

We will prove that $m_k=0$, for all $k$, by proving that for
each $\sigma^{1}_{R}$, all $n_i$ are zero. Fix one $\sigma_R^1$.
Like before, identify $(\sigma_{C}^{1})_i$ with $(i,i+1)$, and
note that
$\Gamma^{(i,i+1)\times\sigma^{1}_{R}}=\Gamma^{\{i+1\}\times\sigma^{1}_
{R}}=\Gamma^{\sigma^{1}_{R}}\cong D_6$, for all $i+1\neq n(\sigma^1_R)$,
because $(\sigma_{C}^{1})_i=(i,i+1)$ is a ray cell having the same
direction as $\sigma^{1}_{R}$. Hence $n_i\in
K_{-1}(\intZ[\Gamma^{\sigma^{1}_{R}}])$ and the boundary maps for
the coefficients are just the identity:

$$\begin{array}{ccc}
K_{-1}(\intZ[\Gamma^{[i,i+1]\times\sigma^{1}_{R}}]) &
\stackrel{\partial}{\longrightarrow} &
 K_{-1}(\intZ[\Gamma^{\{i+1\}\times\sigma^{1}_{R}}]) \\
 \| & & \| \\
 K_{-1}(\intZ[\Gamma^{\sigma^{1}_{R}}]) & \stackrel{id}{\longrightarrow} &
 K_{-1}(\intZ[\Gamma^{\sigma^{1}_{R}}]). \\

\end{array}
$$\\

Since all $\sigma^{2}_{k}$ are special 2-cells, Lemma 5.5
implies that, for $i\geq 0$ and $i+1\neq n(\sigma^{1}_{R})$,
$[i,i+1]\times\overline{\sigma^{1}_{R}}$ and $[i+1,i+2]\times\overline
{\sigma^{1}_{R}}$ are the only closed 2-cells in the family $\{\overline
{\sigma^{2}_{k}}\}$ whose boundaries contain $\{i+1\}\times\overline{\sigma_{R}^{1}}$.
Hence, since we are assuming
$$\partial(\sum_{i}n_i((\sigma^1_C)_i\times\sigma_R^1))=0,$$

\n we have that $n_i=n_{i+1}$. But $[0,1]\times\overline{\sigma^{1}_{R}}$ is the
only closed 2-cell of the family $\{\, [i,i+1]\times\overline{\sigma^{1}_{R}}\,
\}$ whose boundary contains $[\, 0\, ]\times\sigma^{1}_{R}$. Hence
$n_0=0$. Therefore $n_i=0,\, i<n(\sigma^{1}_{R})$. On the other
hand, if $n_{i_0}\neq 0,\, i_0>n(\sigma^{1}_{R})$ we would have
$n_i\neq 0$ for all $i>i_0$, which shows that the sum $\sum_{i}n_i
(((i,i+1)\times\sigma^{1}_{R}))$ is infinite. Hence
$n_i=0$ for all $i$. This proves the lemma.\\

\vspace{1cm}

\centerline{\large{\bf Calculation of the term  ${\cal E}^{2}_{1,-2}$.}}

\vspace{1cm}

Consider the associated cellular
 chain complex, $...\leftarrow
C_{0}\stackrel{\partial}\longleftarrow C_{1}\leftarrow...$, for
${\cal E}^{2}_{1,-2}$, where
 $C_{0}
\hspace{0.2cm}
=\bigoplus_{j}{\tilde
K}_{0}(\intZ[\Gamma^{\sigma_{j}^{0}}])$ and
 $C_{1}=\bigoplus_{i}{\tilde
K}_{0}(\intZ[\Gamma^{\sigma_{i}^{1}}])$.
From Section 4 (isotropy of 1-cells), $\Gamma^{\sigma_{i}^{1}}$ is
trivial or isomorphic to $\intZ_i, D_i, \intZ, D_{\infty},
\intZ_i\times\intZ_{2}, D_i\times\intZ_{2},
\intZ\times\intZ_{2}, D_{\infty}\times\intZ_{2}, i=2,3,4,6$.
By \cite{R} (see also \cite{Ro}), $\tilde K_{0}(\intZ
[F])=0$ for $F$ isomorphic to  $D_i$ or  $\intZ_{i}$,
$i=2,3,4$. By \cite{R}, $\tilde K_{0}(\intZ
[\intZ_6])=0$. By \cite{B2},  $\tilde K_{0}(\intZ
[F])=0$ if $F$ is isomorphic to  $D_{\infty}$ or $\intZ$. By \cite{P}, $\tilde
K_{0}(\intZ [F])=0$ if $F$ is isomorphic to  $D_{\infty}\times\intZ_{2}$,
$\intZ\times\intZ_{2}$, or $D_6$. $\tilde
K_{0}(\intZ [F])$ does not vanish for $F$ isomorphic to
$\intZ_{4}\times\intZ_{2}$, $\intZ_{6}\times\intZ_{2}$,
$D_2\times\intZ_{2}$, $D_{4}\times\intZ_{2}$ or
$D_{6}\times\intZ_{2}$ (see \cite{B2} and \cite{BM}).
Hence it is enough to study the
case $\sigma^1\subset A-Z\cong_{\Gamma}(C^{3}\times \real^{3})-Y$, in
which
$\sigma^1=\sigma^{1}_{C}\times\sigma^{0}_{R}$, where
$\sigma^{1}_{C}$ is a ray cell and $\sigma^{0}_{R}$ is a 0-cell
(in the other cases ${\tilde K}_{0}$ vanishes).\\

In the following lemma we consider $D_{\infty}=\intZ_2*\intZ_2$,
and write $D_{\infty}=\intZ_2^{a}*\intZ_2^{b}$ to distinguish the
two factors. Let $i^a:\intZ_2^{a}\hookrightarrow\intZ_2^{a}*
\intZ_2^{b}$, $i^b:\intZ_2^{b}\hookrightarrow\intZ_2^{a}*
\intZ_2^{b}$ be the canonical inclusions, let $\beta:\intZ_2^{a}*
\intZ_2^{b}\rightarrow\intZ_2^{a}$ be homomorphism such that
$\beta (a)=a$ and $\beta (b)=0$. Then $\beta\circ
i^a=1_{\tiny\intZ_2^{a}}$
and $\beta\circ i^b=0$. Let $F$ be a group.
Define the group homomorphisms \\

\n $\alpha^a:\intZ_2^{a}\times
F\rightarrow(\intZ_2^{a}*\intZ_2^{b})\times F$,\,
$\alpha^b:\intZ_2^{b}\times
F\rightarrow(\intZ_2^{a}*\intZ_2^{b})\times F$ by\\

\n $\alpha^a=i^a\times 1_F$, $\alpha^b=i^b\times 1_F$. Let\\

\n$\alpha^a_{*}:{\tilde K}_0(\intZ[\intZ_2^{a}\times
F])\longrightarrow
\tilde K_0(\intZ[(\intZ_2^{a}*\intZ_2^{b})\times F])$,\\

\n$\alpha^b_{*}:{\tilde K_0}(\intZ[\intZ_2^{b}\times F])
\longrightarrow \tilde K_0(\intZ[(\intZ_2^{a}*\intZ_2^{b})\times
F])$
 denote the induced homomorphisms, at the $\tilde K_0$ level.\\

\n{\bf Lemma 5.7:} \textit{Let $F$ be a  group with $\tilde
K_0(\intZ [F])=0$. Then $\alpha_{*}^{a}$ and $\alpha_{*}^{b}$ are
injective. Moreover
$Im(\alpha_{*}^{a})\cap Im(\alpha_{*}^{b})=\{0\}$.}\\

\n{\bf Proof:} Consider the following  diagram of groups and
homomorphisms

$$\begin{array}{ccccc}
\vspace{0,3cm}
\intZ_2^{a}\times
F & \stackrel{\alpha^{a}}\longrightarrow &
(\intZ_2^{a}*\intZ_2^{b})\times F& \stackrel{\gamma}\longrightarrow &
\intZ_2^{a}\times F.\\
\vspace{0,3cm}
 & & \uparrow{\alpha^{b}} & &  \\
&  &  \intZ_2^{b}\times F &  & \\
\end{array}
$$

\n Here $\gamma=\beta\times 1_F$. Note that
$\gamma\circ\alpha^{a}=1_{{\tiny \intZ_2^{a}}\times F}$, and
$\gamma\circ\alpha^{b}=0$. Applying the $\tilde K_0$ functor, we get

$$\begin{array}{ccccc}
\vspace{0,3cm}
 \tilde K_0(\intZ[\intZ_2^{a}\times
F]) & \stackrel{\alpha_*^{a}}\longrightarrow & \tilde
K_0(\intZ[\intZ_2^{a}*\intZ_2^{b}\times F])&  \stackrel{\gamma_*}
\longrightarrow &
\tilde K_0
(\intZ[\intZ_2^{a}\times F]).\\
\vspace{0,3cm}
 & & \uparrow{\alpha_*^{b}} & &  \\
&  & \tilde K_0(\intZ[\intZ_2^{b}\times F]) &  &\\
\end{array}
$$

\n Then $\gamma_*\circ\alpha_*^{a}$ is the identity. Hence $\alpha_*^{a}$ is injective.
Analogously $\alpha_*^{b}$ is injective. If $x\in
Im(\alpha_{*}^{a})\cap Im(\alpha_{*}^{b})$, then $x=\alpha_*^{a}\,
 y$ and $x=\alpha_*^{b}\, z$. Hence
${\gamma_*}x={\gamma_*}\circ\alpha_*^{a}\, y=y$ and
${\gamma_*}x={\gamma_*}\circ\alpha_*^{b}\, z=0$. Therefore $x=0$. This
proves the lemma.
\\

\n{\bf Lemma 5.8:} \textit{Consider the 1-cells $\sigma_C^1\times\sigma_R^0$,
$\tilde\sigma_C^1\times\tilde\sigma_R^0$ of $A-Z$, where $\sigma_C^1,
\tilde\sigma_C^1 $ are  ray cells. Suppose that
$\overline{(\sigma_C^1\times\sigma_R^0)}\cap\overline{(\tilde\sigma_C^1\times
\tilde\sigma_R^0)}\neq\emptyset$.}\\

\n{\bf (a)} \textit{If $\sigma_R^0=\tilde\sigma_R^0$ then either
$\sigma_C^1$ and $\tilde\sigma_C^1 $ lie in the same ray (i.e.,
determine same direction in $\real^3$) or
$\overline{\sigma_C^1}\cap\overline{\tilde\sigma_C^1}=[0]$}.\\

\n{\bf (b)} \textit{If $\sigma_R^0\neq\tilde\sigma_R^0$ then
$\overline{\sigma_C^1}\cap\overline{\tilde\sigma_C^1}\neq\emptyset$,
$\sigma_C^1$, $\tilde\sigma_C^1$ lie in the same ray and the line
$l$ that contains $\sigma_R^0$ and $\tilde\sigma_R^0$ has the
direction determined by $\sigma_C^1$. Moreover $l\in{\cal L} $ and
$\overline{(\sigma_C^1\times\sigma_R^0)}\cap\overline{(\tilde\sigma_C^1\times
\tilde\sigma_R^0)}=[\, {\bar l}\, ]$.}\\

\n{\bf Proof (a):} If $\sigma_R^0=\tilde\sigma_R^0$, then
\vspace{0,1cm}
$\sigma_C^1$ and $\tilde\sigma_C^1$ lie in the same cone based on
\vspace{0,1cm}
$\sigma_R^0=\tilde\sigma_R^0$. 
Since the collapsing map $\pi : C^n\times\real^n\rightarrow A$ is injective on all cones
$C^n\times\{ v\}$, $v\in\real^3$, we have that
$\overline{(\sigma_C^1\times\sigma_R^0)}\cap\overline{(\tilde\sigma_C^1\times
\tilde\sigma_R^0)}\neq\emptyset$, implies either
$\sigma_C^1\times\sigma_R^0=\tilde\sigma_C^1\times
\tilde\sigma_R^0$ or
$\overline{(\sigma_C^1\times\sigma_R^0)}\cap\overline{(\tilde\sigma_C^1\times
\tilde\sigma_R^0)}=\sigma_C^0\times\sigma_R^0$. If
$\sigma_C^0=[0]$, then
\vspace{0,1cm}
$\overline{\sigma_C^1}\cap\overline{\tilde\sigma_C^1}=[0]$.
Otherwise $\sigma_C^1$ and $\tilde\sigma_C^1$ lie in the same ray,
because
$\overline{\sigma_C^1}\cap\overline{\tilde\sigma_C^1}=\sigma_C^0$ and
$\sigma_C^1$, $\tilde\sigma_C^1$ are ray cells. \\

\n{\bf (b):} If $\sigma_R^0\neq\tilde\sigma_R^0$, consider
the line $l$ that contains $\sigma_R^0$ and $\tilde\sigma_R^0$.
Recall that we collapse only the lines $\bar l$ with $l\in {\cal L}$. Note
that the height at which we collapse depends only on the direction of
the line, and that if a point in a ray is identified with a point
in another ray, then these
\vspace{0,1cm}
rays determine the same direction. Hence if
$\overline{(\sigma_C^1\times\sigma_R^0)}\cap\overline{(\tilde\sigma_C^1\times
\tilde\sigma_R^0)}\neq\emptyset$ with
$\sigma_R^0\neq\tilde\sigma_R^0$, we have that $l\in {\cal L}$.
\vspace{0,1cm}
Moreover
$\overline{(\sigma_C^1\times\sigma_R^0)}\cap\overline{(\tilde\sigma_C^1\times
\tilde\sigma_R^0)}=[\, \bar l\, ]$ and
$\overline{\sigma_C^1}\cap\overline{\tilde\sigma_C^1}\neq\emptyset$.
This proves the lemma.\\

The next lemma implies that ${\cal E}^{2}_{1,-2}=0$.\\

 \n{\bf Lemma 5.9:} \textit{ $\partial_1:C_1\longrightarrow
C_2$ is injective.}\\

\n{\bf Proof:} Let $\sum_{k}m_{k}\sigma^{1}_{k}$ $\in C_1$ with $\partial
( \sum_{k}m_{k}\sigma^{1}_{k} )=0$.
Since $m_{k}\in {\tilde K}_{0}(\intZ[\Gamma^{\sigma_{k}^{1}}])$,
we can suppose that $\sigma_{k}^{1}$ is of the form
$\sigma_{C}^{1}\times\sigma^{0}_{R}$, with $\sigma_{C}^{1}$ a ray
cell (in the other cases ${\tilde
K}_{0}(\intZ[\Gamma^{\sigma_{k}^{1}}])=0$). Thus we can write

$$\sum_{k}m_{k}\sigma^{1}_{k}=\sum_{\sigma^{0}_{R},\,\, r}\left(\sum_{i}n_{i,r}
((\sigma_{C}^{1})_{i,r}\times\sigma^{0}_{R})\right).$$

\n where we assume that $\cup_i{\overline{(\sigma_{C}^{1})}_{i,r}}$ is 
\vspace{0,1cm}
the ray $r$ and that the $(\sigma_{C}^{1})_{i,r}$ are enumerated in such a way
that $(\sigma_{C}^{0})_{i,r}:=\overline{(\sigma_{C}^{1})}_{(i-1),r}\cap\overline
{(\sigma_{C}^{1})}_{i,r}\neq\emptyset$. Note that all
$\Gamma^{(\sigma_{C}^{1})_{i,r}\times\sigma^{0}_{R}}$ are equal, for
all $i$ and $\sigma_{R}^{0}$, $r$ fixed.
We will suppress the subindex $r$ to aliviate the heavy notation, 
e.g we will write $n_{i}$ instead of $n_{i,r}$.\\

We will prove that $m_k=0$, for all $k$, by proving that for each
$\sigma^{0}_{R}$ and $r$ fixed, all $n_i$ are zero.\\

Suppose there is  $\sigma^0_R$ and $r$ for which there is a $n_i$
\vspace{0,1cm}
with $n_i\neq 0$. Since $n_i\in\tilde
K_0(\intZ[\Gamma^{(\sigma^1_C)_i\times\sigma^0_R}])$ then $\tilde
K_0(\intZ[\Gamma^{(\sigma^1_C)_i\times\sigma^0_R}])\neq 0$. Hence
\vspace{0,1cm}
$\Gamma^{(\sigma^1_C)_i\times\sigma^0_R}$ is isomorphic to
$\intZ_2\times\intZ_4, \intZ_2\times\intZ_6, \intZ_2\times D_2,
\intZ_2\times D_4, \intZ_2\times D_6 $, i.e.,
$\Gamma^{(\sigma^1_C)_i\times\sigma^0_R}$ is isomorphic to
$\intZ_2\times F$, where $F$ is isomorphic to $\intZ_4, \intZ_6, D_2, D_4$
or $D_6$.\\

Let $l$ be the line with direction determined by the direction of $r=
\cup_i\overline{(\sigma^1_C)}_i$
and that contains $\sigma^0_R$. Let $\alpha$ be the plane
\vspace{0,1cm}
orthogonal to $l$ that contains $\sigma_{R}^{0}$. Let $G=\Gamma^{(l)}$. Note that
$\Gamma^{(\sigma^1_C)_i\times\sigma^0_R}=G^{\sigma_R^0}$, for all
\vspace{0,1cm}
$i$. Moreover $G^{\sigma_R^0}$ acts by reflection on $l$ (if the
action was trivial,
$G^{\sigma_R^0}=\Gamma^{(\sigma^1_C)_i\times\sigma^0_R}$ would be
\vspace{0,1cm}
isomorphic to $\intZ_i,  D_i, i=2,3,4,6$). Also for
$\Gamma^{(\sigma^1_C)_i\times\sigma^0_R}\cong\intZ_2\times F$, the
$\intZ_2$ factor acts by reflection on $l$, and $F$ acts on the
plane $\alpha$. Note that $F=G^{l}$. Hence we have two cases:
$l/G$ is not compact or $l/G$ is compact.\\

\n{\bf First case:} $l/G$ is not compact.\\

Then $G$ contains no translations. Hence $
G=G^{\sigma_R^0}=\Gamma^{(\sigma^1_C)_i\times\sigma^0_R}\cong\intZ_2\times
F$, which is finite. Therefore, $l\notin {\cal L}$. Since
$r=\cup_i{\overline{(\sigma_{C}^{1})}_i}$ is the ray determined by $l$,
we have
that $(r\times\sigma_R^0)\cap Z=\emptyset$ (we only collapse lines
\vspace{0,1cm}
$\bar l$ with $l\in {\cal L}$). Hence
$\Gamma^{(\sigma^1_C)_i\times\sigma^0_R}=\Gamma^{(\sigma^0_C)_{i+1}\times\sigma^0_R}=G$
and the boundary maps\\
$\tilde
K_0(\intZ[\Gamma^{(\sigma^1_C)_i\times\sigma^0_R}])\rightarrow\tilde
K_0(\intZ[\Gamma^{(\sigma^0_C)_{i+1}\times\sigma^0_R}])$  at the coefficient level
are just the identity.\\

Now recall that we are assuming some $n_i\neq 0\in\tilde K_0(\intZ [G])$.
Since $(r\times \sigma_{R}^{0})\cap Z=\emptyset$, (b) of lemma 5.8
cannot happen. Hence, by (a) of lemma 5.8, we have that
$n_{i+1}\neq 0\in\tilde K_0(\intZ [G]) $. In the same way we
\vspace{0,1cm}
prove that $n_j\neq 0\in\tilde K_0(\intZ [G])$ for all $j\geq i$
which is a contradiction, because the sum above is finite.\\

\n{\bf Second case:} $l/G$ is  compact.\\

Then $G$ contains translations. Moreover, since $G$ acts by
reflection on $l$, we have that the action of $G$ on $l$ is a
dihedral action. Then $G=G^{\sigma_R^0}*_F G^{\tilde\sigma_R^0}$
where $F=G^l$ is isomorphic to $\intZ_4, \intZ_6, D_2, D_4$ or
$D_6$, and $\tilde\sigma_R^0\in l$,
$\tilde\sigma_R^0\neq\sigma_R^0$. Also, since $G$ contains
translations, $l\in {\cal L}$. Thus $l\in {\cal L}_m$ for some $m$.
Since $r=\cup_i\overline{(\sigma^1_C)}_i$ is the ray determined
\vspace{0,1cm}
by $l$, by lemma 5.8 (b) we have that $(r\times\sigma_R^0)\cap
Z=(r\times\tilde\sigma_R^0)\cap Z=(r\times\sigma_R^0)\cap$
\vspace{0,1cm}
$(r\times\tilde\sigma_R^0)=\{[\, \bar
l\, ]\}=[\{m\}\times\sigma_R^0]=$
\vspace{0,1cm}
$[\{m\}\times\tilde\sigma_R^0]$ ($\bar
l$
\vspace{0,1cm}
is collapsed to a point). Write $\sigma^0=\{[\, \bar l\, ]\}$, for
this 0-cell. Note that
\vspace{0,1cm}
$G=\Gamma^{(l)}=\Gamma^{([\{m\}\times\sigma_R^0])}$
and that
$\Gamma^{(\sigma^0_C)_{i}\times\sigma^0_R}=\Gamma^{(\sigma^1_C)_{i+1}
\times\sigma^0_R}=
G^{\sigma^0_R}$
\vspace{0,1cm}
for all $i$ and $i+1\neq m$. Then the boundary  maps (at the coefficient level)
$\tilde K_0(\intZ[\Gamma^{(\sigma^1_C)_{i}\times\sigma^0_R}])\rightarrow\tilde
K_0(\intZ[\Gamma^{(\sigma^0_C)_{i+1}\times\sigma^0_R}])$ are just the
identity for $i+1\neq m$. That is, we have, for $i+1\neq m$, the folowing diagram:

$$\begin{array}{ccc}
\tilde K_{0}(\intZ[\Gamma^{(\sigma^1_C)_i\times\sigma^0_R}]) &
\stackrel{\partial}{\longrightarrow} & \tilde
K_{0}(\intZ[\Gamma^{(\sigma^0_C)_{i+1}\times\sigma^0_R}])
  \\
 \| & & \| \\
 \tilde K_{0}(\intZ [G^{\sigma^0_R}]) & \stackrel{id}{\longrightarrow} &
 \tilde K_{0}(\intZ [G^{\sigma^0_R}]). \\

\end{array}
$$\\
For $i+1=m$, we have

$$\begin{array}{ccc}
\tilde K_{0}(\intZ[\Gamma^{(\sigma^1_C)_{i}\times\sigma^0_R}]) &
\stackrel{\partial}{\longrightarrow} & \tilde K_{0}(\intZ
[\Gamma^{\sigma^0}])
  \\
 \| & & \| \\
 \tilde K_{0}(\intZ [G^{\sigma^0_R}]) &
\stackrel{\alpha_*}{\longrightarrow} &
 \tilde K_{0}(\intZ [G]), \\

\end{array}
$$\\

\n where $\alpha_{*}$ is induced by the inclusion $\alpha:G^{\sigma_r^0}\rightarrow
G=G^{\sigma_R^0}*_F G^{\tilde\sigma_R^0}$.
By Lemma 5.7 $\alpha_*$ is injective. Analogously, $\tilde\alpha_*$
is injective, where $\tilde\alpha:G^{\tilde\sigma_r^0}\rightarrow
G=G^{\sigma_R^0}*_F G^{\tilde\sigma_R^0}$ is the inclusion.\\

\n{\bf Claim:} $Im(\alpha_*)\cap Im(\tilde\alpha_*)=\{0\}$.\\

\n{\bf Proof:} Note that
$|G^{\sigma_R^0}/F|=|G^{\tilde\sigma_R^0}/F|=2$. If
$G^{\sigma_R^0}\cong G^{\tilde\sigma_R^0}\cong\intZ_2\times F$, the
claim follows from Lemma 5.7, because $\tilde K_0(\intZ [F])=0$ for $F$ isomorphic
to $\intZ_4, \intZ_6, D_2, D_4$ or  $D_6$.  In the other cases we
have $\tilde K_0(\intZ [G^{\tilde\sigma_R^0}])=0$. This
proves the claim.\\

Now, recall that we are assuming that some $n_i\neq 0\in \tilde
K_0(\intZ [\Gamma^{(\sigma_C^1)_i\times\sigma_R^0}])=\tilde
K_0(\intZ [G^{\tilde\sigma_R^0}])$. If $i+1\neq m$, Lemma 5.8 (a)
implies that $n_{i}=n_{i+1}$, for $i+1\neq m$. 
In the other case, $i+1=m$, Lemma 5.8 (b) and 
the claim above imply that $n_{m-1}=n_m$. Then $0\neq n_i=n_{i+1}$, for all
$i$, which is a contradiction because the sum above is finite. This
proves the lemma.\\

\section{ Proof of the Results}

\hspace{0.5cm} Here we proof the results mentioned in the introduction.\\

\n{\bf Proof of the Main Theorem.}
First a claim.\\

 \n{\bf Claim} \textit{The terms ${\cal E}^{\infty}_{p,q}$
 vanish if $p+q=-1$ with the exception of $${\cal
E}^{\infty}_{0,-1}=
\bigoplus_{G\in
{VC}_{\infty}(\Gamma)}Wh(G)
.$$}\\

\n{\bf Proof.} To calculate the terms ${\cal E}^{\infty}_{0,-1}$
 we need of the terms  ${\cal
E}^{2}_{2,-2}$ and ${\cal E}^{2}_{3,-3}$. Recall that if
${\sigma^2}$ is a 2-cell in $A$, $\Gamma^{\sigma^2}$ is trivial or isomorphic
to  $\intZ, D_{\infty}, D_{i}, \intZ{_i},
\, i=2,3,4,6$. Then ${\tilde K}_{0}(\intZ [\Gamma^{\sigma^2}])=0$ for
any these groups.
Hence ${\cal E}^{2}_{2,-2}=0$. Recall also that if $\sigma^3$ is a
3-cell in $A$, $\Gamma^{\sigma^3}$ is trivial or isomorphic to
$\intZ_{2}$. Then ${ K}_{-1}(\intZ
[\Gamma^{\sigma^3}])=0$. Hence ${\cal E}^{2}_{3,-3}=0$.
Therefore
$${\cal E}^{2}_{0,-1}={\cal E}^{3}_{0,-1}=...={\cal E}^
{\infty}_{0,-1}=
\bigoplus_{G\in
{VC}_{\infty}(\Gamma)}Wh(G)
.$$
\n In Section 5 we proved that ${\cal E}^{2}_{1,-2}$ and ${\cal
E}^{2}_{2,-3}$ vanish. Hence ${\cal E}^{\infty}_{1,-2}$ and ${\cal
E}^{\infty}_{2,-3}$ vanish. This
proves the claim.\\

Recall that our spectral
sequence, \vspace{0.1cm} for $p+q=-1$,  converges to
\vspace{0,1cm}
$H_{-1}(A/\Gamma,{\cal P}_{*}(\rho)):=\pi_{-1}({\mathbb
H}(A/\Gamma,{\cal P}_{*}(\rho)))$, where ${\cal P_*}( )$ is  the
stable pseudoisotopy  functor (see \cite{FJ3}). 
By the claim, propositions 1.2, 1.3, 1.4 we have that $$Wh(\Gamma)=
Wh(\pi_{1}(X))=\pi_{-1}({\cal P}_{*}(X))= \bigoplus_{G\in
{VC}_{\infty}(\Gamma)}Wh(G),$$ where $X$ is a space such that $\pi_{1}(X)=\Gamma$.\\

It remains to prove that the sum in the formula is finite.
This is implied by the following lemma:\\

\n{\bf Lema 6.1:} \textit{ Let $\Gamma$ be a 3-crystallographic group.
Then, the
subset $\{G<\Gamma; G\in VC_{\infty}, Wh(G)\neq 0 \}$
is finite.}\\

\n{\bf Proof:} If $G\in VC_{\infty}$, with $Wh(G)\neq 0$, then $G$
leaves invariant only one  line in $\real^3$ (see Lemma 5.3). Let
$l_G$ be the
line invariant by $G$.\\

\n{\bf Claim 1:} \textit{ $l_G$ is contained in the 1-skeleton of
the
triangulation ${\cal T}_R$ of $\real^3$.}\\

\n{\bf Proof of  Claim 1:} Since $G\in VC_{\infty}$, we have the
sequence
$0\rightarrow F\rightarrow G\rightarrow
H\rightarrow 0$, with $H$ isomorphic to $\intZ$ or $D_{\infty}$ and $F$ a
 finite subgroup of  $\Gamma$. Recall that  $F$ fixes $l_{G}$
pointwise. Recall also that the $\Gamma$-action on $\real^3$ is simplicial.
Moreover, we can assume (subdividing if necessary) that $\gamma\sigma =\sigma$ implies
$\gamma|_{\sigma}=1|_{\sigma}$, where $\sigma$ is a simplex and $\gamma\in\Gamma$.
Hence:\\

If $\sigma^3$ is a open 3-cell  in ${\cal T}_R$ and
$l_G\cap\sigma^3\neq\emptyset$, then $\gamma\sigma^3=\sigma^3$ for all
 $\gamma\in F$. Hence $F$ fixes pointwise $\real^3$ and follows that $F$
is the trivial group. Therefore, $G$ is isomorphic to   $\intZ$ or
$D_{\infty}$ and
$Wh(G)=0$.\\

If $\sigma^2$ is a open 2-cell in ${\cal T}_R$ and $l_G\cap\sigma^2\neq\emptyset$,
we have that $\gamma\sigma^2=\sigma^2$ for all
$\gamma\in F$. Then $F$ fixes  pointwise a plane $\alpha$
in $\real^3$ and acts (at most) by reflection in $\alpha^{\perp}$.
Therefore $F$ is the  trivial group or isomorphic to  $\intZ_2$.
By propositions 3.2 and 3.3, we have that $G$ is isomorphic to one of
the following groups: $\intZ\times\intZ_2$,
$D_2*_{\tiny\intZ_2}D_2$, $D_2*_{\tiny\intZ_2}\intZ_4$,
$\intZ_4*_{\tiny\intZ_2}\intZ_4$. Consider the group
$G=G_1*_{\tiny\intZ_2}G_2$, where $G_i$, $i=1,2$ are finite. By
\cite{Wa}, we obtain the exact sequence
$$Wh(\intZ_2)\rightarrow Wh(G_1)\oplus Wh(G_2)\rightarrow
Wh(G)\rightarrow \tilde K_0(\intZ[\intZ_2]).$$ Since
$Wh(\intZ_2)$, $Wh(\intZ_4)$, $Wh(D_2)$ and $\tilde
K_0(\intZ[\intZ_2])$ are trivial, we have that $Wh(D_2*_{\tiny\intZ_2}D_2)$,
$Wh(D_2*_{\tiny\intZ_2}\intZ_4)$,
$Wh(\intZ_4*_{\tiny\intZ_2}\intZ_4)$ are trivial.  Also, $Wh(\intZ\times\intZ_2)$ is trivial
(see \cite{P}). This proves the claim.\\

\n{\bf Claim 2:} \textit{Let $\Gamma^{(l_{1})}$ and \,
$\Gamma^{(l_{2})}$ be
the isotropy groups  of the lines $l_{1}$ and $l_{2}$
respectively. Assume that  $l_i$ is the only line left invariant by $\Gamma^{(l_{i})}$ $i=1,2.$ 
Then there is $\gamma\in\Gamma$ such that $\gamma l_{1}=l_{2}$ if and only
if  $\Gamma^{(l_{1})}$ and
$\Gamma^{(l_{2})}$ are  conjugate.}\\

\n{\bf Proof of  Claim 2:} If $\gamma\in\Gamma$ is such that
$\gamma\Gamma^{(l_{1})}\gamma^{-1}=\Gamma^{(l_{2})}$ and
$g_{1}\in\Gamma^{(l_{1})}$, $g_{2}\in\Gamma^{(l_{2})}$ with
$g_{2}=\gamma g_{1}\gamma^{-1}$, we have that: $g_{2}(\gamma
l_{1})=(\gamma g_{1}\gamma^{-1})(\gamma l_{1})=\gamma
g_{1}l_{1}=\gamma l_{1}$. Therefore $\Gamma^{(l_{2})}$ fixes the line
$\gamma l_{1}$. Since $\Gamma^{(l_{2})}$ fixes only one line,
follows that $\gamma l_{1}=l_{2}$.  Conversely, if $\gamma l_1 =l_2$ and
$g_{1}\in\Gamma^{(l_{1})}$, then $\gamma
g_{1}\gamma^{-1}l_{2}=\gamma g_{1}l_{1}=\gamma l_{1}=l_{2}$. Hence
$\Gamma^{(l_{1})}$ is conjugate
to $\Gamma^{(l_{2})}$. This proves the claim.\\

Recall that the  $\Gamma$-action on $\real^3$ is cocompact. Then,
there is a finite subcomplex $D\subset\real^3$, such that for all 
$x\in \real^3$ there are $x_D\in D$ and $\gamma\in \Gamma$ with $\gamma x_D=x$.
 Hence, if $l$ is a line in $\real^3$, there are a line $l_D$ and
$\gamma\in \Gamma$ with $l_D\cap D\neq\emptyset$ and
$\gamma l_D=l$. Since $D$ is finite, $D$ intercepts only a
finite number of 1-cells of the triangulation ${\cal T}_R$. By
Claims 1 and 2, follows that, modulo conjugation, there is a finite number of
groups $G\in VC_{\infty}$
with $Wh(G)\neq 0$. This proves the lemma and the Main Theorem.\\

\n{\bf Proof of Corollary 1:} Follows directly from the Main Theorem.\\

The following result was used in the example mentioned in the introduction.\\

\n{\bf Lemma 6.2:} \textit{
$Nil_1(\intZ[\,  \intZ_2\times\intZ_2\, ])$ is infinitely generated.}\\

\n{\bf Proof:} Consider the Cartesian square with all maps
surjective:

$$\begin{array}{ccc}
\intZ[G] & {\longrightarrow} & \intZ[\intZ_2]\\
 {\downarrow} & & {\downarrow} \\
 \intZ[\intZ_2] &{\longrightarrow} & F_2[\intZ_2]. \\
\end{array}
$$

\n where $F_2$ denotes the field with two elements. By \cite{Hr},
$Nil_1(\intZ[\intZ_2])=0$. Hence the Mayer-Vietoris sequence of
that square produce an epimorphism
$$Nil_2(F_2[\intZ_2])\longrightarrow Nil_1(\intZ [G])\longrightarrow 0.$$
Since $Nil_2(F_2[\intZ_2])$ is non trivial (this follows from
 van der Kallen technique \cite{VK}), we have that  $ Nil_1(\intZ [G])\neq 0$.
Since  $ Nil_1(\intZ [G])\neq 0$, we have that $ Nil_1(\intZ [G])$ is
infinitely generated (see \cite{F}). This proves the lemma.\\

To prove  corollary 2 we need two lemmata.\\

\n{\bf Lemma 6.3:} \textit{Let $\Gamma$ be a n-crystallographic
group and $G$ a maximal virtually infinite cyclic group of
$\Gamma\times\intZ$. Then (at least) one of the following holds:}\\

\n\textit{{\bf (1)} $G$ is isomorphic to a virtually
infinite cyclic subgroup of $\Gamma$.}\\

\n\textit{{\bf (2)} $G=F\times\intZ$, $F<\Gamma$, with $F$ maximal finite in $\Gamma$.}\\

\n{\bf Proof:}  Let
$\pi:\Gamma\times\intZ\rightarrow\Gamma$ be the projection. Then
$Ker(\pi)=\intZ$. We have two possibilities: $\pi(G)$ finite or
$\pi(G)$ infinite. \\

\n{\bf Claim:} If $\pi(G)$ is infinite then $\pi(G)\cong G$.\\

\n{\bf Proof of  Claim:} In fact, if $\pi|_{G}$
is not one-to-one then we get a exact sequence $0\rightarrow
Ker(\pi|_{G})\rightarrow G\rightarrow\pi(G)\rightarrow 0$ with
$Ker(\pi|_{G})$ and $\pi(G)$
infinite. This is impossible since $G$ is virtually  cyclic. This
proves the claim.\\

Suppose $\pi(G)$ finite. Write $F=\pi(G)$. Then $G<F\times\intZ$.
But $G$ is maximal, hence $G=F\times\intZ$. Certainly, F has to be
maximal. This proves the
lemma.\\

\n{\bf Lemma 6.4:} \textit{  
$F$ is maximal finite in $\Gamma$ if and only if $F\times\intZ$
is maximal  virtually infinite cyclic  in $\Gamma\times\intZ$.}\\

\n{\bf Proof:} If $F\times \intZ$ is maximal  virtually infinite cyclic 
in $\Gamma\times\intZ$ then, certainly, $F$ is maximal finite in $\Gamma$.\\
Suppose now that $F$ is maximal finite in $\Gamma$ and let $H$ be a virtually
infinite cyclic  subgroup of $\Gamma$, with $F\times \intZ \subset H$.
Let $\pi: \Gamma\times\intZ\rightarrow \Gamma$ be the projection.
If $\pi(H)$ is infinite then we have the exact sequence
$$0\rightarrow (Ker \pi)\cap H\rightarrow H\rightarrow \pi(H)\rightarrow 0$$
This is a contradiction because $\{ 0\}\times\intZ =(Ker\pi )\cap H$ and $H$ is
virtually infinite cyclic. Follows that $\pi (H)$ is finite.
Since $F$ is maximal finite in $\Gamma$, we have that $F=\pi (H)$.
Therefore $H=F\times\intZ$. This proves the lemma.\\

We will use the following fact  (see \cite{P}):\\

\n{\bf (*)} $Wh(G)={\tilde K}_{0}(G)=0$, for $G$ a finite or virtually infinite cyclic subgroup of
2-crystallographic group.\\

Recall that $F(\Gamma) $ denotes the set of conjugacy classes of maximal
finite subgroups of $\Gamma$.\\

\n{\bf Proof of  corollary 2:} Let $\Gamma$ be a
2-crystallographic group.  Bass-Heller-Swan-formula (see [28]
p.152) implies that
\vspace{0,1cm}
$Wh(\Gamma\times\intZ)=Wh(\Gamma)\oplus{\tilde
K}_{0}(\intZ[\Gamma])\oplus 2Nil_{1}(\intZ[\Gamma])$. But
$Wh(\Gamma),{\tilde K}_{0}(\intZ[\Gamma])$ vanish for any
\vspace{0,1cm}
2-crystallographic group (see (*) above).  Hence
$Wh(\Gamma\times\intZ)=2Nil_{1}(\intZ[\Gamma])$. 
Note that $F$ is conjugate to $F'$ in $\Gamma$ if and only if
$F\times\intZ$ is conjugate to $F'\times\intZ$ in $\Gamma\times\intZ$.
Hence, by the Main Theorem, Lemmata  6.3, 6.4, and the Bass-Heller-Swan formula 
we have
$$Wh(\Gamma\times\intZ)=
\bigoplus_{G\in
{VC}_{\infty}(\Gamma\times{\tiny\intZ})}Wh(G)=$$
$$\bigoplus_{F\in
F(\Gamma)}Wh(F\times\intZ)\\
 =\bigoplus_{F\in
F(\Gamma)}2Nil_{1}(\intZ[F])=2\, (\bigoplus_{F\in
F(\Gamma)}Nil_{1}(\intZ[F])\,\, ).$$

\pagebreak

A. Alves

Departamento de Matematica

Universidade Federal de Pernambuco

Cidade Universitaria, Recife PE 50670-901

Brazil

e-mail: almir@dmat.ufpe.br
\vspace{.3in}

P. Ontaneda

Departamento de Matematica

Universidade Federal de Pernambuco

Cidade Universitaria, Recife PE 50670-901

Brazil

e-mail: ontaneda@dmat.ufpe.br or pedro@math.binghamton.edu

\end{document}